\documentclass[10pt,reqno]{amsart}
\usepackage{amssymb,color}
\usepackage{amsmath,amssymb,graphicx}
\usepackage{accents}
\usepackage{bbm}
\usepackage{lipsum}
\usepackage{array}
\usepackage{multirow}

\usepackage{tabularx,ragged2e,booktabs,caption}
\usepackage{here}
\usepackage[toc,page]{appendix}

\usepackage[colorlinks=true, citecolor=blue, linkcolor=blue]{hyperref}
\usepackage{amsthm}
\definecolor{c20}{rgb}{0.,0.7,0.}
\definecolor{c30}{rgb}{0.,0.,1.}
\definecolor{c40}{rgb}{1,0.1,0.7}
\definecolor{c50}{rgb}{1,0,0}
\definecolor{c60}{rgb}{0,0.9,0.1}

\newcommand{\kb}[1]{\boldsymbol{#1}}
\newcommand{\vk}[1]{\kb{#1}}

\newcommand{\BQN}{\begin{eqnarray}}
\newcommand{\EQN}{\end{eqnarray}}
\newcommand{\BQNY}{\begin{eqnarray*}}
\newcommand{\EQNY}{\end{eqnarray*}}

\def\green#1{\textcolor{c30}{#1}}
\def\green#1{{#1}}

\newcommand{\BS}{\begin{sat}}
\newcommand{\ES}{\end{sat}}
\newcommand{\BT}{\begin{theo}}
\newcommand{\ET}{\end{theo}}
\newcommand{\BK}{\begin{korr}}
\newcommand{\EK}{\end{korr}}

\newcommand{\BD}{\begin{de}}
\newcommand{\ED}{\end{de}}

\newcommand{\BIT}{\begin{itemize}}
\newcommand{\EIT}{\end{itemize}}
\newcommand{\BDI}{\begin{description}}
\newcommand{\EDI}{\end{description}}

\newcommand{\BRM}{\begin{remark}}
\newcommand{\ERM}{\end{remark}}

\newcommand{\BEL}{\begin{lem}}
\newcommand{\EEL}{\end{lem}}

\newtheorem{theo}{Theorem}[section]
\newtheorem{sat}[theo]{Proposition}
\newtheorem{de}[theo]{Definition}
\newtheorem{lem}[theo]{Lemma}

\newtheorem{example}[theo]{Example}
\newtheorem{korr}[theo]{Corollary}
\newtheorem{remark}[theo]{Remark}

\newcommand{\prooftheo}[1]{ \textbf{Proof of Theorem} \ref{#1} }
\newcommand{\proofprop}[1]{\textbf{Proof of Proposition} \ref{#1}}

\newcommand{\COM}[1]{}

\newcommand{\QED}{\hfill $\Box$ \\}

\date{}

\def\x{\vk{x}}

\begin{document}

\title{On some multivariate Sarmanov mixed Erlang reinsurance risks: aggregation and capital allocation}

\author{Gildas Ratovomirija}
\address{Gildas Ratovomirija, Department of Actuarial Science, University of Lausanne, UNIL-Dorigny 1015 Lausanne, Switzerland 
and 	Vaudoise Assurances, Place de Milan CP 120, 1001 Lausanne, Switzerland}

\author{Maissa Tamraz}
\address{Maissa Tamraz, Department of Actuarial Science, University of Lausanne, UNIL-Dorigny 1015 Lausanne, Switzerland}

\author{Raluca Vernic}
\address{Raluca Vernic, Faculty of Mathematics and Computer Science, Ovidius University of Constanta, 124 Mamaia Blvd., 900527 Constanta, Romania
and Institute for Mathematical Statistics and Applied Mathematics,
Calea 13 Septembrie 13, 050711 Bucharest, Romania}

\bigskip

\date{\today}
 \maketitle

\bigskip
{\bf Abstract:} 
Following some recent works on risk aggregation and capital allocation for mixed Erlang risks joined by Sarmanov's multivariate distribution, in this paper we present some  closed-form formulas for the same topic by considering, however, a different kernel function for Sarmanov's distribution, not previously studied in this context. The risk aggregation and capital allocation formulas are derived and numerically illustrated in the general framework of stop-loss reinsurance, and then in the particular case with no stop-loss reinsurance. A discussion of the dependency structure of the considered distribution, based on Pearson's correlation coefficient, is also presented for different kernel functions and illustrated in the bivariate case.\\


					\section{Introduction}
Modern risk management usually involves complex dependent risk factors. In this respect,  several regulations were put in place in order to assess the minimum capital requirement, namely the Economic Capital (EC) that insurance and reinsurance companies are constrained to hold according to their risk exposures. In practice, the EC is evaluated  by means of risk measures on the aggregated risk, so that the companies will be covered from unexpected large losses. 

For instance, the EC under the Solvency II framework for EU countries focuses on a Value-at-Risk (VaR) approach at a tolerance level of 99.5\% of the aggregated risk over a one year period, while in Switzerland, the EC under the Swiss Solvency Test (SST) is based on the Tail-Value-at-Risk (TVaR) approach at a 99\% confidence level of the aggregated risk over a one year period. 
Since the EC quantified in the latter reflects the aggregate capital needed to cover the entire loss of a company, it is also of interest to study how this capital should be allocated among the different risk factors (e.g., lines of business) in the insurance and reinsurance companies, in other words, how much amount of capital each individual risk contributes to the aggregated EC. This allows the risk managers to identify and monitor conveniently their risks.
An extensive literature has been developed on capital allocation techniques from which we shall restrict to the TVaR method (see  \cite{dhaene2012optimal}, \cite{tasche2004allocating}  and the references therein for an overview of the existing methods). Our choice is motivated by the fact that Artzner \cite{artzner1999coherent} discussed the properties of the  VaR risk measure and showed that it fails to fulfill all the axioms of a coherent risk measure (hence, it might not be a reasonable tool for capital allocation), while the TVaR fulfills all the axioms and, moreover, provides information on the tail of the distribution.\\

Therefore, the main task of actuaries is to choose an appropriate model for the multivariate risk factors, namely the dependence structure model and the  distributions of the marginals. 
The aim of this contribution is to  address  risk aggregation and TVaR capital allocation for insurance and reinsurance mixed Erlang risks whose dependency is governed by the Sarmanov distribution with a certain expression of the kernel functions. This study comes along the lines of some recent contributions:  Vernic \cite{vernic2015capital} considered capital allocation based on the TVaR rule for the Sarmanov distribution with exponential marginals; Cossette et al. \cite{Cossette_al13} used the Farlie-Gumbel-Morgenstern (FGM) distribution to model the dependency between mixed Erlang distributed risks and applied it to capital allocation and risk aggregation;  Hashorva and Ratovomirija \cite{Hashorva_Rija14} and Ratovomirija \cite{Ratovomirija2016} presented aggregation and capital allocation in insurance and reinsurance for mixed Erlang distributed risks joined by the Sarmanov distribution with a specific kernel function different from the one considered in this study. Note that the choice of the Sarmanov and mixed Erlang distributions is not incidental, these distributions gained a lot of interest in the actuarial literature lately: for the Sarmanov distribution, see e.g.,  \cite{Yang_al13}, \cite{hernandez_al12}, \cite{boucher2015sarmanov}, \cite{vernic2016Sarmanov}, while for the mixed Erlang distribution we refer to  \cite{Lee_Lin10}, \cite{Willmot_al10}, \cite{Lee_al12} or \cite{Willmot_and_woo_14}.  One key advantage of the Sarmanov distribution is its flexibility to join different types of marginals and its allowance to obtain exact results. An interesting property of the mixed Erlang distributions is the fact that many risk related quantities, such as TVaR, have an analytical form. \\

This paper is organized as follows: in the second section, we present some preliminaries on the Sarmanov distribution, on the TVaR capital allocation problem and on the mixed Erlang distribution, supplemented with several lemmas on this last distribution that will be needed for the proofs of the main results.
Section 3 contains the main results on risk aggregation and capital allocation for the stop-loss reinsurance, which are also particularized in the case without stop-loss reinsurance; the main formulas of this section are illustrated with some numerical examples. The paper ends with two appendices: the first one discusses and compares the dependence structure of the bivariate Sarmanov distribution with mixed Erlang marginals and different kernel functions, providing upper and lower bounds for the corresponding Pearson correlation coefficient, while the second appendix contains all the proofs of the theoretical results.

					\section{Preliminaries}
\subsection{Multivariate Sarmanov distribution}
The Sarmanov distribution caught the interest of many researchers in different fields. It was first introduced by Sarmanov \cite{Sarmanov66} in the bivariate case, then extended by Lee \cite{Lee96} to the multivariate case. Its applications in many insurance contexts show its flexible structure  when modeling the dependence between multivariate risks given the distribution of the marginals. For instance, Abdallah et al. \cite{boucher2015sarmanov} used a bivariate Sarmanov distribution to model the dependence within or between  lines of business through calendar years, accident years and development years in the loss reserving framework, while Hernandez et al. \cite{hernandez_al12} developed a new Sarmanov family with beta and gamma marginals used for the computation of the  Bayes premium in a collective risk model. \\

According to \cite{Sarmanov66}, the joint probability density function (pdf) of a bivariate Sarmanov distribution is defined as follows
\BQN \label{eq: densityBivariate}
	h(x_1,x_2) =f_1(x_1) f_2(x_2) (1+\alpha_{1,2}  \phi_{1}(x_{1})\phi_{2}(x_{2})),x_1,x_2\in \mathbb{R},
\EQN

where for $i=1,2,f_i$ are the densities of the marginals and $\phi_i$ are kernel functions assumed to be bounded, non-constant and satisfying the following conditions 
\BQNY
	\mathbb{E}(\phi_{i}(X_{i} ))=0, i=1,2, \quad 1+\alpha_{1,2} \phi_{1}(x_{1})\phi_{2}(x_{2}) \geq 0,\forall x_1,x_2\in \mathbb{R}.
\EQNY
Lee \cite{Lee96} introduced general methods for the choice of $\phi$. Yang and Hashorva \cite{Yang_al13} considered the case where $\phi$ depends on some function $g$, being expressed as follows
$$ \phi(x)=g(x)- \mathbb{E}(g(X)),  \text{  where  } \mathbb{E}(g(X)) < \infty.$$
In the context of risk aggregation and capital allocation, Hashorva and Ratovomirija \cite{Hashorva_Rija14} assumed that $g(x)= e^{-x}$, Vernic \cite{vernic2015capital} studied the case where the marginals are exponentially distributed, while Cossette et al. \cite{Cossette_al13} used the 
FGM distribution with mixed Erlang marginals (the FGM is a special case of the Sarmanov distribution for $g(x)=2(1-F(x)),$ with $F$ denoting the distribution function of the marginal). Thus, in the sequel, we consider the following kernel function
\BQN\label{eq:condkern1}
\phi_{i}(x_{i})= f_i(x_i) - \mathbb{E}(f_i(X_i)),
\EQN
in which case the range of $\alpha_{1,2}$ is given by  
\BQN \label{eq:Alpha12Bound}
	  \frac{-1}{\max \{ \gamma_{1}  \gamma_{2},(M_1-\gamma_{1})(M_2- \gamma_{2}) \}}
	  \leqslant  \alpha_{1,2}  \leqslant
	   \frac{1}{\max \{ \gamma_{1} (M_2- \gamma_{2}),(M_1-\gamma_{1}) \gamma_{2} \}},
\EQN
where $\gamma_{i}= \mathbb{E}(f_i(X_i))$ and $M_i= \underset{x \in \mathbb{R}}{\text{max }}f_i(x)$,$ i=1,2.$ Moreover, we shall work with a generalization of the above distribution to the multivariate case,  see \cite{Lee96}. In this respect, for simplicity, we denote, in the rest of the paper, by $\vk{X}:=(X_1,\ldots,X_n)$ an \textit{n}-variate random vector, by $\vk{x}:=(x_1,\ldots,x_n)$ an \textit{n}-dimensional vector (e.g., the observations on $\vk{X}$) and we let $I_{n}=1,\ldots,n$. Therefore, we shall model the dependency between the risks $X_i$ having pdf $f_i,i \in I_{n},$ via the multivariate Sarmanov distribution having the following pdf
\BQN \label{eq:Sarmanovdensity}
h(\x)=\prod _{i=1}^n f_i(x_i)\left( 1+\sum_{1 \leq j<l \leq n} \alpha_{j,l}  \phi_{j}(x_{j})\phi_{l}(x_{l})\right), \x \in \mathbb{R}^n,
\EQN

where $\phi_i$  are the non-constant kernel functions defined in (\ref{eq:condkern1}) 
and $ \alpha_{j,l}$  are real numbers satisfying the condition
\BQN \label {eq: condkern2}
1+\sum_{1 \leq j<l \leq n} \alpha_{j,l}  \phi_{j}(x_{j})\phi_{l}(x_{l}) \ge 0.
\EQN

\BRM
It should be noted that a more general expression of the Sarmanov density for the multivariate case can be written as follows
\BQN \label{eq: Sarmanov_general}
h(\vk{x})= \prod_{i=1}^{n} f_i(x_i) \left( 1+ \sum_{l=2}^{n}\sum_{1 \leq j_1<j_2<\ldots<j_l\leq n} \alpha_{j_1,\ldots,j_l} \prod_{k=1}^l\phi_{j_k}(x_{j_k})\right),\x \in \mathbb{R}^n, 
\EQN 
such that  $\mathbb{E}(\phi_{i}(X_{i} ))=0$ 
and $ 1+ \sum_{l=2}^{n}\sum_{1 \leq j_1<j_2<\ldots<j_l\leq n} \alpha_{j_1,\ldots,j_l} \prod_{k=1}^l\phi_{j_k}(x_{j_k}) \geq 0.$
However, (\ref{eq: Sarmanov_general}) requires the estimation of all the dependence parameters, which is in general very complex. Thus, it is often assumed that $ \alpha_{j_1,\ldots,j_l}=0$ for $l \geq 3$, see \cite{mari2001correlation}. For simplicity, in this paper, we consider the Sarmanov density defined in (\ref{eq:Sarmanovdensity}). \\
\ERM

\subsection{Mixed Erlang distributions}
The mixed Erlang distribution has many attractive distributional properties when modeling the claim sizes of an insurance portfolio, see, e.g., \cite{Willmot_al10}, and the dependence between multivariate insurance risks, see \cite{Lee_al12}. Actually, during these past few years, modeling the dependence of multivariate mixed Erlang risks raised the interest of many researchers. For instance, Cossette et al. \cite{Cossette_al13} modeled the dependence of multivariate mixed Erlang risks using the FGM distribution and looked at its applications in the risk management framework. Moreover, Hashorva and Ratovomirija \cite{Hashorva_Rija14} and Ratovomirija \cite{Ratovomirija2016}
studied the dependence of mixed Erlang risks governed by the Sarmanov distribution in the context of capital allocation and risk aggregation in insurance and reinsurance.\\

In this regard, we define the pdf of a mixed Erlang distribution denoted $ME(\beta,\underline{Q})$ by
\BQN
f(x,\beta,\underline{Q})=\sum_{k=1}^{\infty} q_k w_k(x,\beta),x\geq0 ,  
\EQN
where $w_k(x,\beta)= \frac{\beta^k x^{k-1} e^{-\beta x}}{(k-1)!}$ is the pdf of an Erlang distribution with $\beta>0$  the scale parameter, $k \in \mathbb{N}^*$ the shape parameter and $\underline{Q}=(q_1,q_2,\ldots)$ is a vector of non-negative mixing probabilities such that $\sum_{k=1}^{\infty} q_k=1$.
We denote by $W_k$ the distribution function (df) of the Erlang distribution and by $\overline{W}_k$ its corresponding survival (tail) function given, respectively, by
\BQN\notag
W_k(x,\beta)=1-e^{-\beta x}\sum_{j=0}^{k-1}\frac{(\beta x)^j }{j!}, \quad \overline{W}_k(x,\beta)=e^{-\beta x}\sum_{j=0}^{k-1}\frac{(\beta x)^j }{j!},x\geq0 .
\EQN
Thus, the mixed Erlang df can be expressed in terms of the Erlang df as follows
\BQN
F(x,\beta, \underline{Q})= \sum_{k=1}^{\infty} q_k W_k(x,\beta)
= 1- e^{-\beta x} \sum_{k=1}^{\infty} q_k \sum_{j=0}^{k-1} \frac{(\beta x)^j}{j!},x\geq0 .
\EQN
Moreover, the expected value of this distribution is $\mu= \frac{1}{\beta} \sum_{k=1}^{\infty} kq_k $.

In addition, we present some distributional properties and useful results for the mixed Erlang distributions.  

\begin{lem}\label{lem:tranformation ME}
Let $X \sim ME(\beta,\underline{Q})$ with pdf $f(x, \beta, \underline{Q}) $ and $\mathbb{E} (f(X,\beta, \underline{Q})) <\infty$. Then $c(x, \beta, \underline{Q}): =\frac{f(x, \beta, \underline{Q})^2}{\mathbb{E} \left( f(X,\beta, \underline{Q})\right)} $ is again a pdf of a mixed Erlang distribution with mixing probabilities $ V(\underline{Q})=(v_1,v_2,\ldots) $  and scale parameter $2\beta$, i.e., we have 
\BQNY 
c(x, \beta,\underline{Q})= \sum_{k=1}^{\infty} v_k w_k(x, 2\beta)=f\left( x, 2\beta,V( \underline{Q} )\right)  ,
\EQNY
where 
\BQN \label{eq:vk}
v_k&=& \frac{\sum_{i=1}^{k}\begin{pmatrix} {k-1} \\ {i-1} \end{pmatrix}\frac{q_i  q_{k+1-i}}{2^{k} }}{ \sum_{i=1}^{\infty}  \sum_{j=1}^{\infty}   \begin{pmatrix} {i+j-2} \\ {i-1} \end{pmatrix}\frac{q_i  q_j}{2^{i+j-1} }}.
\EQN 
\end {lem}

The proof of the above lemma is given in the Appendix \ref{SectionProofs}. We shall use the notation $ \tilde{\mu} $ for the expected value corresponding to the pdf $c(\cdot, \beta, \underline{Q})$ defined in this lemma, i.e., 
\BQN \label{mu_i}
\tilde{\mu}=\frac{1}{2\beta}\sum_{k=1}^{\infty} kv_k.
\EQN 
The following results have already been developed in \cite{Cossette_al13} and \cite{Ratovomirija2016}.

\begin{lem}\label{lem:tranformation ME g=x}
Let $X \sim ME(\beta,\underline{Q})$ with pdf $f(x,\beta,\underline{Q})$. Then $f^{G}(x,\beta,\underline{Q}):=\frac{x f(x,\beta,\underline{Q})}{\mathbb{E}(X)}$ is equal to the pdf $f(x,\beta,G(\underline{Q}))$ of a mixed Erlang distribution with mixing probabilities $G(\underline{Q})=(g_1,g_2,\ldots)$ given by

\BQN\notag
g_{k}=
  \begin{cases}
 0 & \text{if }  k=1\\
\frac{(k-1)q_{k-1}}{\sum_{j=1}^{\infty}jq_j} &\text{if }  k=2,3,\ldots
\end{cases} .
\EQN \\
\end {lem}

\begin{lem}\label{lem:tranformation ME with different scale parameters}
Let $X \sim ME(\beta_1,\underline{Q})$. Then it follows that for any  positive constant $\beta_2$  such that $\beta_2 \geq \beta_1$ , we have $X \sim ME(\beta_2,\Psi(\underline{Q})),$
where the elements of $\Psi(\underline{Q})=(\psi_1,\psi_2,\ldots)$ are given by
\BQN\notag
\psi_k=\sum_{i=1}^k q_i  \begin{pmatrix} {k-1} \\ {k-i} \end{pmatrix} \left( \frac{\beta_1}{\beta_2}\right) ^{i}\left( 1-\frac{\beta_1}{\beta_2}\right) ^{k-i}, k\geq 1.
\EQN
\end {lem}

\begin{lem}\label{lem:Aggregation}
Let $X_1,X_2$  be two independent mixed Erlang random variables (r.v.s) such that $X_i \sim ME(\beta,\underline{Q}_i) ,i=1,2$. Then $S_2:=X_1+X_2 \sim ME\left( \beta, \underline{\Pi} \left( \underline{Q}_1, \underline{Q}_2 \right) \right) $ with the components of $ \underline{\Pi} \left( \underline{Q}_1, \underline{Q}_2 \right) $ given by
	\BQNY
		\pi_l \left( \underline{Q}_1, \underline{Q}_2 \right) =
			\left\{
			 	\begin{array}{rcl}
         			0  & \mbox{for} & l=1 \\
         			\sum_ {j=1}^{l-1} q_{1,j}  q_{2,l-j}& 		    \mbox{for} & l>1
              	\end{array}
        	\right. .
	\EQNY
\BRM \label{rem:aggregME}
According to Remark 2.1 in \cite{Cossette_al12}, the result in Lemma \ref{lem:Aggregation} can be extended to $S_n:= \sum_{i=1}^n X_i$, given that $X_1,\ldots,X_n$ are independent r.v.s and $X_i\sim ME(\beta,\underline{Q}_i) $   for $i\in I_{n}$. Thus, $S_n \sim ME\left( \beta, \underline{\Pi} \left( \underline{Q}_1,\ldots, \underline{Q}_n \right) \right) $, where the  mixing weights are determined iteratively as follows
\BQNY
		\pi_l \left( \underline{Q}_1,\ldots, \underline{Q}_{n+1}\right) =
			\left\{
			 	\begin{array}{rcl}
         			0  & \mbox{for} & l=1,\ldots,n \\
         			\sum_ {j=n}^{l-1} \pi_j\left( \underline{Q}_1,\ldots, \underline{Q}_{n}\right) q_{n+1,l-j}& 		    \mbox{for} & l=n+1,n+2,\ldots
              	\end{array}
        	\right. .
	\EQNY
\ERM
\end{lem}	

\begin{lem}\label{lem:Reinsriskdf}
Given $d>0$ and the r.v. $X \sim ME(\beta,\underline{Q})$ , the df of $Y := (X-d)_{+}$ can be expressed as 
\BQN
F_Y(y)=F_{X}\left( y+d\right) =F_{X}\left( d\right)
+H_{X}\left( y,d\right) ,y\geq 0,  \label{FY}
\EQN
where 
$$H_{X}\left( y,d \right) :=\mathbb{P} \left( 0<Y\leq y\right) =\sum_{k=0}^{\infty
}\Delta _{k}\left( d,\beta ,\text{\underline{Q}}\right) W_{k+1}\left( y,\beta
\right) ,$$
with
$$\Delta _{k}\left( d,\beta ,\text{\underline{Q}}\right) =\beta
^{-1}\sum_{j=0}^{\infty }q_{j+k+1}w_{j+1}\left( d,\beta \right) .
$$
Moreover, defining $ U_X(y,d):=\int_y^ \infty u \frac{\partial}{\partial u}H_{X}\left( u,d \right) du, $ it also holds that
\BQNY	
U_X(y,d)= \frac{1}{\beta}
	 \sum_ {k=0}^{\infty} 
			 (k+1)   \Delta_k(d,\beta,\underline{Q})
			 \overline{W}_{k+2}(y,\beta), \green{ y >0}.
\EQNY		
\end {lem}

The following result is proved in Section \ref{SectionProofs}. We introduce the convention that an empty product equals 1.

\begin{lem} \label{lem:Aggreinsrisk}

Consider the independent r.v.s $X_{i}\sim ME\left( \beta ,\underline{Q}%
_{i}\right) $, let $d_{i}>0$ and $Y_{i}=\left( X_{i}-d_{i}\right) _{+},i\in
I_{n}.$ Then the df of $R_{n}=\sum_{i=1}^{n}Y_{i}$ can be written as
\BQN
F_{R_{n}}\left( y\right) =\prod\limits_{i=1}^{n}F_{X_{i}}\left( d_{i}\right)
+\sum_{k=1}^{n}\sum_{1\leq j_{1}<...<j_{k}\leq
n}H_{X_{j_{1}}+...+X_{j_{k}}}\left( y,d_{j_{1}},...,d_{j_{k}}\right)
\prod\limits_{i\in I_{n}\setminus \left\{ j_{1},...,j_{k}\right\}
}F_{X_{i}}\left( d_{i}\right) ,y\geq 0,  \label{FRn}
\EQN
where, for $k\geq 1,$ 
\BQNY
H_{\sum_{i=1}^{k}X_{i}}\left( y,d_{1},...,d_{k}\right)
&:=& \mathbb{P}\left(
\bigcap\limits_{i=1}^{k}\left( X_{i}>d_{i}\right) ,\sum_{i=1}^{k}\left(
X_{i}-d_{i}\right) \leq y\right)  \notag \\
&=& \sum_{h_{1}=0}^{\infty
}...\sum_{h_{k}=0}^{\infty }\Delta _{h_{1}}\left( d_{1},\beta ,\underline{Q}%
_{1}\right) \cdot ...\cdot \Delta _{h_{k}}\left( d_{k},\beta ,\underline{Q}%
_{k}\right) W_{\sum_{i=1}^{k}h_i+k}\left( y,\beta \right) .
\EQNY
Moreover, if $U_{\sum_{i=1}^{k}X_{i},X_{k+1}}(y,d_1,\ldots,d_{k+1}):=\int_y^ \infty \int_0^s
 u \frac{\partial}{\partial u} H_{X_{k+1}}(u,d_{k+1}) \left[ \frac{\partial}{\partial v}
 H_{\sum_{i=1}^{k}X_{i}}\left( v,d_{1},...,d_{k}\right)\right]  _{v=s-u} du ds, $ then
\BQNY
	U_{\sum_{i=1}^{k}X_{i},X_{k+1}}(y,d_1,\ldots,d_{k+1})=
	\frac{1}{\beta}
\sum_{h_1=0}^{\infty}... \sum_{h_{k+1}=0}^{\infty} 
\left( h_{k+1} +1 \right)  
\Delta_{h_1}(d_1,\beta,\underline{Q}_1) \cdot...\cdot \Delta_{h_{k+1}}(d_{k+1},\beta,\underline{Q}_{k+1}) 
\overline{W}_{\sum_{i=1}^{k+1}h_i+k+2}(y,\beta).
\EQNY
\end {lem}

\subsection{TVaR capital allocation}
As mentioned in the introduction, it is of great interest for insurance and reinsurance companies to quantify the total capital required for the safety of the company, and also to determine the part of this capital to be allocated to each risk/portfolio in order to cover its loss. Among the capital allocation techniques discussed in the literature, we shall consider the TVaR rule. In order to present the allocation formulas, we recall the definitions of the VaR and TVaR risk measures for a risk $X$ and tolerance level $ p \in \left( 0,1\right) $, i.e.,
\BQNY
VaR_{p}\left( X\right) =\min
\left\{ x\left\vert F_{X}\left( x\right) \geq p \right. \right\},TVaR_{p }\left( X\right) =\mathbb{E}\left( X\left\vert
X>VaR_{p }\left( X\right) \right. \right). 
\EQNY
Let $X_{i} $ denote the $i$th risk r.v. of an insurance portfolio and let $S=\sum_{i=1}^{n}X_{i}$ represent the aggregate risk of the portfolio. Then, if the total risk capital is evaluated as $TVaR_{p}\left( S\right)$, the TVaR capital allocation rule naturally allocates to the $i$th risk 
\BQNY
C_{i}(p)=TVaR_{p}\left( X_{i}, S\right) :=\mathbb{E}\left(
X_{i}\left\vert S>VaR_{p }\left( S\right) \right. \right) ,
\EQNY
which can be rewritten as
\BQN \label{eq:ReinsTVar}
C_{i}(p)=\frac{1}{1-p}\mathbb{E}\left( X_{i} \mathbbm{1}_{\left\lbrace  S>VaR_{p}(S)\right\rbrace }\right),
\EQN
where $ \mathbbm{1}_A $ denotes the indicator function of the set $A$. Clearly, $TVaR_{p}(S)=\sum_{i=1}^n C_{i}(p).$

\section{Main results }

\subsection{Joint distribution of aggregate Sarmanov risks}

We consider $n$ insurance portfolios where each portfolio consists of $k_1,\ldots, k_n$ risks, respectively. We denote by  $S_{i}=\sum_{j=1}^{k_i}X_j^{(i)}$ the aggregate risk of portfolio $ i $, where $X_j^{(i)}$ is the $j$th  individual risk from the $i$th portfolio having pdf $f_j^{(i)}, j=1,\ldots,k_i,i \in I_n$. We assume that the joint distribution of $\vk{X} :=\left( X_1^{(1)},\ldots, X_{k_1}^{(1)}; \ldots;X_{1}^{(n)},\ldots, X_{k_n}^{(n)} \right) $ is governed  by Sarmanov's distribution with the pdf as defined in \eqref{eq:Sarmanovdensity} and fulfilling  \eqref{eq:condkern1} and \eqref{eq: condkern2} for  the kernel  functions $\phi$, i.e., in this case,  
\BQN \label{eq: Sarmanovrv}
h(\vk{x})&=& \prod_{i=1}^{n} \prod_{j=1}^{k_i} f_j^{(i)}\left( x_j^{(i)}\right)  \left[  1+ \sum_{1\leq a< b\leq n}
\sum_{s=1}^{k_a} \sum_{t=1}^{k_b} \alpha_{s,t}^{(a,b)} \phi_{s}^{(a)} \left( x_{s}^{(a)}\right) \phi_{t}^{(b)} \left( x_{t}^{(b)}\right) \right. \notag \\
&&\left. + \sum_{a=1}^{n} \sum_{1\leq s< t\leq k_a}\alpha_{s,t}^{(a)} \phi_{s}^{(a)} \left( x_{s}^{(a)}\right) \phi_{t}^{(a)} \left( x_{t}^{(a)}\right) \right] , 
 \EQN 
where $ \vk{x}= \left( x_1^{(1)},\ldots, x_{k_1}^{(1)}, \ldots,x_{1}^{(n)},\ldots, x_{k_n}^{(n)} \right) $. 

Next, we present the joint density of $\vk{S} = (S_{1}, \ldots, S_{n} ) $ under these assumptions.
 
  \begin{theo} \label{theorem:JointVectors}
The joint pdf of $\vk{S}$ can be expressed as follows 
\BQNY
f_{\vk{S}}(s_1,\ldots,s_n)=\prod_{i=1}^{n}f_{S_{i}}(s_{i})+\sum_{1\leq a\leq b\leq
n}\sum_{s=1}^{k_{a}}\sum_{t=T\left( s,a,b\right) }^{k_{b}}\alpha
_{s,t}^{(a,b)}\prod_{i=1}^{n}\tilde{f}_{S_{i;s,t}^{\left( a,b\right) }}(s_{i}),
\EQNY
where $T\left( s,a,b\right) =\max \left\{ 1,\left( s+1) \mathbbm{1}_{\left( a=b\right)
}\right) \right\},\alpha_{s,t}^{(a,a)}=\alpha_{s,t}^{(a)}$, 
\BQN\label{eq:convol}
f_{S_{i}}= f_1^{(i)}* \ldots* f^{(i)}_{k_i}, \tilde{f}_{S_{i;s,t}^{\left( a,b\right) }} =\tilde{f}_{1;s,t}^{(i;a,b)}* \ldots* \tilde{f}_{k_i;s,t}^{(i;a,b)},i\in I_n,
\EQN

and, for  $i\in I_n, j=1,\ldots,{k_i},$ 
$$\tilde{f}_{j;s,t}^{(i;a,b)}\left( x\right) =\left\{ 
\begin{array}{ll}
\left( \phi _{j}^{(i)}f_{j}^{(i)}\right) \left( x\right)  & if\ (i,j)\in
\left\{ \left( a,s\right) ,\left( b,t\right) \right\}  \\ 
f_{j}^{(i)}\left( x\right)  & otherwise%
\end{array}%
\right. .$$

\end{theo}
\BRM
It should be noted that  Ratovomirija \cite{Ratovomirija2016} provided a general expression for the joint density of $\vk{S}$ in the particular case when $k_1=\ldots=k_n=k$. 
\ERM

Next, we derive a special case of Theorem \ref{theorem:JointVectors} where we assume that all the marginals are mixed Erlang distributed, i.e., 
 $X_j ^{(i)} \sim ME ( \beta_j^{(i)}, \underline{Q}_j^{(i)})$ with $\underline{Q}_j^{(i)}=(q_{j,1} ^{(i)}, q_{j,2}^{(i)}, \ldots),j=1,\ldots,k_{i},i\in I_n$.
Moreover,   individual risks within and across the portfolios are considered to be joined by Sarmanov's distribution with the joint pdf specified in \eqref{eq: Sarmanovrv} and kernel functions 
 $\phi_j^{(i)}(x)= f^{(i)}_j(x)- \mathbb{E}\left( f^{(i)}_j\left( X_j^{(i)}\right) \right) .$
We denote $$\vk{X} =\left( X_1^{(1)},\ldots, X_{k_1}^{(1)}; \ldots;X_{1}^{(n)},\ldots, X_{k_n}^{(n)} \right) 
     \sim SME_{\zeta}(\vk{\beta}, \vk{\underline{Q}},\vk{\alpha}),$$  
where $\zeta=\sum_{i=1}^n k_i,\vk{\beta}=\Bigl(\beta_1^{(1)},\ldots, \beta_{k_1}^{(1)};\ldots; \beta_{1}^{(n)},\ldots, \beta_{k_n}^{(n)} \Bigr)$, $\vk{\underline{Q}}=
 \Bigl(\underline{Q}_1^{(1)},\ldots, \underline{Q}_{k_1}^{(1)};\ldots; \underline{Q}_{1}^{(n)},\ldots, \underline{Q}_{k_n}^{(n)} \Bigr)$ and $ \vk{\alpha} $ consists of all the $\alpha-$coefficients of the Sarmanov pdf \eqref{eq: Sarmanovrv}. In the following, for simplicity, we also denote $\gamma_{j}^{(i)} = \mathbb{E}\left( f^{(i)}_j\left( X_j^{(i)}\right) \right)$ assuming it exists. \\
  
 \begin{sat}\label{proposition:survival}
If $\vk{X} \sim SME_{\zeta}(\vk{\beta}, \vk{\underline{Q}},\vk{\alpha})$
      with $\beta_{k_n}^{(n)} \ge \beta_j^{(i)},$ for $j=1,\ldots, k_i,i\in I_n,$  then the df of $\vk{S}$ is given by
\BQNY
		F_{\vk{S}}\left( \vk{s}\right)= \xi_n
\prod\limits_{j=1}^{n}F_{S_{j}^{\left( 1 \right) }}\left(
s_{j}\right) -\sum_{1\leq a\leq b\leq n}\sum_{s=1}^{k_{a}}\sum_{t=T\left(
s,a,b\right) }^{k_{b}}\alpha_{s,t}^{(a,b)}\gamma _{s}^{\left( a\right) }\gamma
_{t}^{\left( b\right) } \left( \prod\limits_{j=1}^{n}F_{S_{j;s}^{\left( 2;a\right)
}}\left( s_{j}\right) +\prod\limits_{j=1}^{n}F_{S_{j;t}^{\left(
2;b\right) }}\left( s_{j}\right)
-\prod\limits_{j=1}^{n}F_{S_{j;s,t}^{\left(3;a,b\right) }}\left(
s_{j}\right) \right) ,
	\EQNY

where
\BQN \label{eq:csi}	
\xi_n= 1+\sum_{1\leq a\leq b\leq
n}\sum_{s=1}^{k_{a}}\sum_{t=T\left( s,a,b\right) }^{k_{b}}\alpha_{s,t}^{(a,b)}\gamma
_{s}^{\left( a\right) }\gamma _{t}^{\left( b\right) },
\EQN
\BQN \label{eq:defS1}	
	 S^{(1)}_{j} 
	 \sim ME \Bigl(2\beta_{k_n}^{(n)}, \underline{\Pi} \left( \Psi(\underline{Q}_{1} ^{(j)}), \ldots,  \Psi(\underline{Q}_{k_j} ^{(j)})  \right)  \Bigr) , 
\EQN
\BQN \label{eq:defS3}	
	 S_{j;s,t}^{\left(3;a,b\right) } 
	 \sim ME \left( 2\beta_{k_n}^{(n)}, \underline{\Pi} \left( \Psi \left( M_{s,t}^{(a,b)} (
	 \underline{Q}_{1} ^{(j)}) \right), \ldots,\Psi \left( M_{s,t}^{(a,b)} (
	 \underline{Q}_{k_j} ^{(j)}) \right)  \right) \right)  , 
\EQN	
\BQN \label{eq:defS2}
   S_{j;s}^{\left( 2;a\right)}=S_{j;s,s}^{\left(3;a,a\right) },
\EQN
and 
\BQN \label{eq:def Q}
M_{s,t}^{(a,b)} \left( \underline{Q}_{l} ^{(j)} \right) =
 		\left\{
			 \begin{array}{lcl}
			 V\left(  \underline{Q}_{l} ^{(j)} \right)  & \mbox{if} & (j,l)\in \left\lbrace (a,s),(b,t)\right\rbrace  \notag \\
         		 \underline{Q}_{l} ^{(j)} & \mbox{otherwise} &          			 
              \end{array}              	
      \right. .
\EQN

The components of $ V$ are defined in Lemma \ref{lem:tranformation ME}, the elements of 
 $ \Psi$ in Lemma \ref{lem:tranformation ME with different scale parameters} and the ones of $ \underline{\Pi}$ are given in Remark \ref{rem:aggregME}.
 \end{sat}

\subsection{Stop-loss mixed Erlang reinsurance risks with Sarmanov dependence}
In this section, we study the effect of mixed Erlang distributed risks on reinsurance. In order to mitigate  their risks, insurers enter into reinsurance agreements. There are several types of reinsurance contracts. However, we shall only consider the stop-loss reinsurance. In a stop-loss reinsurance contract, the reinsurer pays the part of the loss that is greater than a certain positive amount $d$ (the deductible). In the following, we shall provide the distribution of the aggregated loss of several reinsurance portfolios in the stop-loss framework, and determine the amount of  capital to be allocated to each reinsurance portfolio under the TVaR allocation principle. \\
In this respect, we consider $n$ insurance portfolios as defined in the last section with aggregated losses $(S_{1},\ldots,S_{n})$ subject to the deductibles $\vk{d}=(d_1,\ldots,d_n)$ on the reinsured amounts $(T_{1},\ldots,T_{n})$, where the $T_{i}$'s, $ i \in I_n $, are defined as follows 
\BQNY
\begin{aligned}
T_{i} &=& (S_{i}-d_i)_+ 
&=&
\begin{cases}
0 & \text{ if } S_{i} \leq d_i\\
S_{i} - d_i & \text{ if } S_{i} > d_i \\
\end{cases}
\end{aligned}.
\EQNY \\
Hereafter, we shall denote by $R_{n}=  \sum_{i=1} ^n T_{i} $ the aggregated reinsurance stop-loss risk.   

\BS\label{prop3}
Let $(X_1^{(1)},\ldots, X_{k_1}^{(1)};\ldots; X_{1}^{(n)},\ldots, X_{k_n}^{(n)} )
     \sim SME_{\zeta}(\vk{\beta}, \vk{\underline{Q}},\vk{\alpha})$
      with $\gamma_{j}^{(i)} < \infty$ and $\beta_{k_n}^{(n)} \ge \beta_j^{(i)},j=1,\ldots, k_i, $
       and let $d_i >0$ for   $i \in I_{n}$. Then the df of $R_{n}$ is given by

\BQN 
F_{R_{n}}\left( y\right)  &=&F_{\mathbf{S}}\left( \mathbf{d}\right) +\sum_{k=1}^{n}\sum_{1\leq j_{1}<...<j_{k}\leq n} \left[ \xi_n H_{S_{j_{1}}^{\left( 1\right) }+...+S_{j_{k}}^{\left( 1\right) }}\left(
y,d_{j_{1}},...,d_{j_{k}}\right) \prod\limits_{i\in I_{n}\setminus \left\{ j_{1},...,j_{k}\right\}
}F_{S_{i}^{\left( 1\right) }}\left( d_{i}\right)\right.  \notag \\
&&- \sum_{1\leq a \leq b\leq
n}\sum_{s=1}^{k_{a}}\sum_{t=T\left( s,a,b\right) }^{k_{b}}\alpha_{s,t}^{(a,b)}
\gamma _{s}^{\left( a\right) }\gamma _{t}^{\left( b\right) } 
\left( H_{S_{j_{1};s}^{\left( 2;a\right) }+...+S_{j_{k};s}^{\left( 2;a\right)
}}\left( y,d_{j_{1}},...,d_{j_{k}}\right) \prod\limits_{i\in I_{n}\setminus
\left\{ j_{1},...,j_{k}\right\} }F_{S_{i;s}^{\left( 2;a\right) }}\left(
d_{i}\right) \right. \notag \\
&&+H_{S_{j_{1};t}^{\left(
2;b\right) }+...+S_{j_{k};t}^{\left( 2;b\right) }}\left(
y,d_{j_{1}},...,d_{j_{k}}\right) \prod\limits_{i\in I_{n}\setminus \left\{
j_{1},...,j_{k}\right\} }F_{S_{i;t}^{\left( 2;b\right) }}\left( d_{i}\right) \notag\\
&&\left. \left. -H_{S_{j_{1};s,t}^{\left( 3;a,b\right) }+...+S_{j_{k};s,t}^{\left(
3;a,b\right) }}\left( y,d_{j_{1}},...,d_{j_{k}}\right) \prod\limits_{i\in
I_{n}\setminus \left\{ j_{1},...,j_{k}\right\} }F_{S_{i;s,t}^{\left(
3;a,b\right) }}\left( d_{i}\right) \right) \right],y\geq 0  ,\label{eq:FRnProp}
\EQN
with $ H $ defined in Lemmas \ref{lem:Reinsriskdf}- \ref{lem:Aggreinsrisk}.
\ES

Next, we shall consider capital allocation under the TVaR principle for the reinsurance risks corresponding to the $n$ portfolios defined above. Let $C_{i}(p)$ be the amount of capital to be allocated to portfolio $i,i \in I_n$, as defined in \eqref{eq:ReinsTVar}. The following result holds.

\BS\label{prop4} 
Let $(X_1^{(1)},\ldots, X_{k_1}^{(1)};\ldots; X_{1}^{(n)},\ldots, X_{k_n}^{(n)} )
     \sim SME_{\zeta}(\vk{\beta}, \vk{\underline{Q}},\vk{\alpha})$ such that  $\gamma_{j}^{(i)} < \infty$ and  $\beta_{k_n}^{(n)} \ge \beta_j^{(i)},$ for $i\in I_n$, $j=1,\ldots, k_i $.
       Let $d_i >0,i\in I_n$ and set $x_p:= VaR_p(R_n)$. Then the capital allocated to portfolio $l$ under the TVaR rule is
\BQNY
C_{l}(p) &=& \frac{1}{1-p} \sum_{k=0}^{n-1}\sum_{\substack{ 1\leq j_{1}<...<j_{k}\leq n \\ l\notin
\left\{ j_{1},...,j_{k}\right\} }}\left[ \xi
_{n}U_{\sum_{i=1}^{k}S_{j_{i}}^{\left( 1\right) },S_{l}^{\left( 1\right)
}}\left( x_{p},d_{j_{1}},..,d_{j_{k}},d_{l}\right) \prod\limits_{i\in
I_{n}\backslash \left\{ j_{1},...,j_{k},l\right\} }F_{S_{i}^{\left( 1\right)
}}\left( d_{i}\right) \right.  \\
&&-\sum_{1\leq a\leq b\leq n}\sum_{s=1}^{k_{a}}\sum_{t=T\left( s,a,b\right)
}^{k_{b}}\alpha _{s,t}^{\left( a,b\right) }\gamma _{s}^{\left( a\right)
}\gamma _{t}^{\left( b\right) }\left( U_{\sum_{i=1}^{k}S_{j_{i};s}^{\left(
2;a\right) },S_{l;s}^{\left( 2;a\right) }}\left(
x_{p},d_{j_{1}},..,d_{j_{k}},d_{l}\right) \prod\limits_{i\in I_{n}\backslash
\left\{ j_{1},...,j_{k},l\right\} }F_{S_{i;s}^{\left( 2;a\right) }}\left(
d_{i}\right) \right.  \\
&&+U_{\sum_{i=1}^{k}S_{j_{i};t}^{\left( 2;b\right) },S_{l;t}^{\left(
2;b\right) }}\left( x_{p},d_{j_{1}},..,d_{j_{k}},d_{l}\right)
\prod\limits_{i\in I_{n}\backslash \left\{ j_{1},...,j_{k},l\right\}
}F_{S_{i;t}^{\left( 2;b\right) }}\left( d_{i}\right)  \\
&&\left. \left. -U_{\sum_{i=1}^{k}S_{j_{i};s,t}^{\left( 3;a,b\right)
},S_{l;s,t}^{\left( 3;a,b\right) }}\left(
x_p,d_{j_{1}},..,d_{j_{k}},d_{l}\right) \prod\limits_{i\in I_{n}\backslash
\left\{ j_{1},...,j_{k},l\right\} }F_{S_{i;s,t}^{\left( 3;a,b\right)
}}\left( d_{i}\right) \right) \right] ,
\EQNY
where, by convention, when $k=0,$ we consider only one term in the sum 
$\sum_{\substack{ 1\leq j_{1}<...<j_{k}\leq n \\ l\notin \left\{ j_{1},...,j_{k}\right\} }} $ in which each component of the type $U_{\sum_{i=1}^{k}S_{j_{i}},S_{l}}$ is replaced with $U_{S_l}$.
\ES

\begin{example}\label{example3.6}
Let $S_{1}$ and $S_{2}$ be the aggregate risks of two insurance portfolios consisting of $k_1=2$ and $k_2=3$ mixed Erlang distributed risks, respectively, with $\beta_3^{(2)} > \beta_j^{(i)}, j=1,\ldots,k_i $ and $i=1,2$. Hence,  $S_{1}=X_1^{(1)}+X_2^{(1)}$ and $S_{2}=X_1^{(2)}+X_2^{(2)}+X_3^{(2)}.$
Following Propositions \ref{proposition:survival}-\ref{prop3}, the distribution of the aggregate stop-loss reinsurance risk $R_2= T_{1} + T_{2}  $, where $ T_{i}=(S_{i} - d_i)_+ , i=1,2$,  is given by 
\BQNY
F_{R_{2}}( y)  &=& 
 \xi_2 \Bigl(F_{S_{1}^{\left( 1\right) }}\left( d_{1}\right)F_{S_{2}^{\left( 1\right) }}\left( d_{2}\right)+ H_{S_{1}^{( 1) }}(y,d_{1}) F_{S_{2}^{\left( 1\right) }}\left( d_{2}\right)+H_{S_{2}^{\left( 1\right) }}\left( y,d_{2}\right) F_{S_{1}^{\left( 1\right) }}\left( d_{1}\right)+ H_{S_{1}^{\left( 1\right) }+S_{2}^{\left( 1\right) }}\left(y,d_{1},d_{2}\right)\Bigr) \\
&&- \sum_{1\leq a\leq b\leq 2}\sum_{s=1}^{k_{a}}\sum_{t=T( s,a,b) }^{k_{b}}\alpha_{s,t}^{(a,b)}
\gamma _{s}^{\left( a\right) }\gamma _{t}^{\left( b\right) } 
\left[  \sum_{(i,j) \in \{(a,s),(b,t)\}} \left( F_{S_{1;j}^{\left( 2;i\right) }}\left(d_{1}\right) F_{S_{2;j}^{\left( 2;i\right) }}\left(d_{2}\right) + H_{S_{{1};j}^{\left( 2;i \right) }}\left( y,d_{{1}}\right) F_{S_{2;j}^{\left( 2;i\right) }}\left(d_{2}\right) \right. \right. \notag \\
&&\left. + H_{S_{{2};j}^{\left( 2;i\right) }}\left( y,d_{{2}}\right) F_{S_{1;j}^{\left( 2;i \right) }}\left(d_{1}\right)+H_{S_{{1};j}^{\left( 2;i\right) }+S_{{2};j}^{\left( 2;i\right) }}\left( y,d_{{1}},d_{{2}}\right)\right)  - F_{S_{1;s;t}^{\left( 3;a;b \right) }}\left(d_{1}\right)F_{S_{2;s;t}^{\left( 3;a;b \right) }}\left(d_{2}\right) \notag \\
&&\left.  -H_{S_{{1};s,t}^{\left( 3;a,b\right) }}\left( y,d_{{1}}\right) F_{S_{2;s,t}^{\left(3;a,b\right) }}\left( d_{2}\right)- H_{S_{{2};s,t}^{\left( 3;a,b\right) }}\left( y,d_{{2}}\right) F_{S_{1;s,t}^{\left(3;a,b\right) }}\left( d_{1}\right)- H_{S_{{1};s,t}^{\left( 3;a,b\right) }+S_{{2};s,t}^{\left( 3;a,b\right) }}\left( y,d_{1},d_{2}\right)\right] .\notag\\
\EQNY

Furthermore, if $TVaR_{p}(R_2)$ is the total risk capital needed to cover $R_2$, in light of Proposition \ref{prop4}, the contribution of $T_{i}$ to this capital is expressed as follows
  \BQNY
 	C_{i}(p) &=& \frac{1}{1-p}\left\lbrace   \xi_2 \left( U_{S_i^{(1)}}(x_p,d_i)F_{S_j^{(1)}}(d_j)+ U_{S_j^{(1)},S_i^{(1)}}(x_p,d_j,d_i) \right) -\sum_{1\leq a\leq b \leq 2}\sum_{s=1}^{k_a} \sum_{t=T(s,a,b)}^{k_b} \alpha_{s,t}^{(a,b)} \gamma_s^{(a)} \gamma_t^{(b)} \right.\\
&&\times \left[ \sum_{(k,l) \in \{(a,s),(b,t)\}} \left( U_{S_{i;l}^{(2,k)}}(x_p,d_i)
F_{S_{j;l}^{(2;k)}}(d_j)+ U_{S_{j;l}^{(2;k)},S_{i;l}^{(2,k)}}(x_p,d_j,d_i)\right) \right. \\
&&\left. \left. -U_{S_{i;s,t}^{(3;a,b)}}(x_p,d_i)F_{S_{j;s,t}^{(3;a,b)}}(d_j)-U_{S_{j;s,t}^{(3;a,b)},S_{i;s,t}^{(3;a,b)}}(x_p,d_j,d_i)\right] \right\rbrace ,i\neq j \in \{1,2\} ,
  \EQNY
 where $ \xi_2$ is defined in (\ref{eq:csi}) for $n=2$, while $ H $ and $U$ are given in Lemmas \ref{lem:Reinsriskdf}- \ref{lem:Aggreinsrisk}. \\

\textbf{Numerical illustration.} To numerically illustrate the just mentioned formulas, in Table \ref{table:quantmeasI} we present the parameters of the individual risks  $X_j^{(i)}$ of the two portfolios, where $j=1,\ldots,k_i$ and $ i=1,2$, together with some related statistical measures (for simplicity, only two decimal places were retained). 
\begin{center}
	\begin{tabular}{|c|c|c|c|c|c|c|c|}
\hline
& $X_j^{(i)}$ & $\beta_j^{(i)}$& $\underline{Q}_j^{(i)}$& Mean & Variance & Skewness & Kurtosis \\		
\hline
	Portfolio I & $X_1^{(1)}$ & 0.12& (0.4,0.6) &13.33&127.78&1.55 &6.50\\
                                &$X_2^{(1)}$ & 0.14& (0.3,0.7) &12.14&97.45 & 1.49& 6.28\\
\hline
	Portfolio II&$ X_1^{(2)}$ &0.15 & (0.5,0.5) &10.00 &77.78&1.62 &6.80\\
           & $X_2^{(2)}$ &0.16& (0.8,0.2) &7.50 & 53.13& 1.88&8.16 \\
           &$ X_3^{(2)} $ &0.18 & (0.55,0.45) &8.06 &52.39& 1.66&6.97\\
\hline
		\end{tabular}
	\captionof{table}{Statistical measures for the individual risks  $X_j^{(i)}, j=1,\ldots,k_i,i=1,2$ (Example \ref{example3.6}).} \label{table:quantmeasI}
\end{center}
Moreover, we assume that the Sarmanov parameters $\alpha_{i,j}$ are as follows
\begin{center}
	\begin{tabular}{c c c c c }
 $ \alpha_{1,2}^{(1)}$= 16,  & $ \alpha_{1,1}^{(1,2)}$=  8, & $ \alpha_{1,2}^{(1,2)}$=5, 
 & $ \alpha_{1,3}^{(1,2)}$=  2,&  $ \alpha_{2,1}^{(1,2)}$=8 , \\
 $ \alpha_{2,2}^{(1,2)}$= 5 ,  &   $ \alpha_{2,3}^{(1,2)}$= 2, & $ \alpha_{1,2}^{(2)}$= 15,
 & $ \alpha_{1,3}^{(2)}$= 17 ,  &   $ \alpha_{2,3}^{(2)}$= 16. \\  
\end{tabular}
\end{center}
Under the stop-loss reinsurance framework, we considered the values $d_1=50$ and $d_2=45$ for the deductibles of Portfolios I and II, respectively. Table \ref{table:reinsallcapital} describes the allocated capitals $C_i(p),i=1,2,$ required to cover the losses of both portfolios after application of the deductibles, as well as the capital needed to cover the loss $R_2$ of the whole reinsured portfolio. We considered several values for the tolerance level $p$.
\begin{center}
	\begin{tabular}{|c||c||c|c|c|}
\hline
p(\%) & $VaR_{p}(R_2)$ & $C_1=TVaR_{p}(T_1,R_2)$ & $C_2= TVaR_{p}(T_2,R_2)$& $ TVaR_{p}(R_2)$\\		
\hline
90.00 &5.03 & 7.33 &8.37 &15.70\\	
\hline
92.50 &8.24& 8.85 &9.90  &18.75 \\	
\hline
95.00 &12.65&11.07  &11.90  &22.97 \\		
\hline
97.50 & 19.96&15.08  &14.96  &30.04\\			
\hline
99.00 &29.31 &20.82  &18.34  &39.16\\		
\hline
99.90 &51.88 & 37.35 &  24.05&61.40 \\	
\hline	
		\end{tabular}
	\captionof{table}{Capital allocated to Portfolios I and II (Example \ref{example3.6}).} \label{table:reinsallcapital}
\end{center}
Table \ref{table:reinsallcapital} shows that for a tolerance level $p \geq 97.5\%$, Portfolio I is riskier than Portfolio II as more capital is needed to cover the losses (this can be explained by the fact that both risks in Portfolio I has higher expected values and variances than the risks in Portfolio II); however, for $p \leq 95\%$, more capital is allocated to Portfolio II. 
\end{example}

\subsection{Particular case: mixed Erlang risks with Sarmanov dependency}
We shall now consider the same setting as before, but in the particular case with only one insurance portfolio, no reinsurance and no deductible. For simplicity, we denote by $\vk{X}=(X_1, \ldots, X_k)$ the $k$ individual risks with joint distribution governed by the $ k$-variate Sarmanov distribution with kernel functions $\phi_j(x_j)=f_j(x_j)-\gamma_j,$ where $\gamma_j=\mathbb{E}(f_j(X_j)),X_j \sim ME(\beta_j,\underline{Q}_j),j \in I_k,$ and we denote by $S=\sum_{j=1}^k X_j $ the aggregate risk of the portfolio. Hence, $\vk{X} \sim SME_k ( \vk{\beta} ,\vk{\underline{Q}},\vk{\alpha}),$ where $  \vk{\beta} =(\beta_1,\ldots,\beta_k),\ \vk{\underline{Q}}=(\underline{Q}_1,\ldots,\underline{Q}_k)$ and $\vk{\alpha}=(\alpha_{i,j})_{1 \leq i<j \leq k}$. Next, we are going to present the distribution of the aggregate risk $S$ that can easily be derived from Proposition \ref{proposition:survival}.

\BS\label{prop1}
 Let $\vk{X} \sim SME_k (\vk{\beta},\vk{\underline{Q}},\vk{\alpha}),$ where $\beta_j \leq \beta_k$ for $j=1,\ldots,k-1$. Then the distribution of the aggregate risk $S$ is given by
\BQNY
 F_{S} (u)&=&\xi_k F_{S^{\left( 1 \right) }}\left(u\right) 
-\sum_{1\leq s < t\leq k} \alpha_{s,t} \gamma_{s}\gamma_{t}
 \left( F_{S_{s}^{\left( 2 \right)}}\left( u \right) +F_{S_{t}^{\left(2\right) }}\left( u\right)
-F_{S_{s,t}^{\left(3\right) }}\left(u\right) \right) ,
\EQNY

where $\xi_k= 1+\sum_{1\leq s < t\leq k} \alpha_{s,t} \gamma_{s}\gamma_{t},$ while
\BQNY 	
	 S^{(1)}  \sim ME \Bigl(2\beta_{k}, \underline{\Pi} \left( \Psi(\underline{Q}_{1}), \ldots,  \Psi(\underline{Q}_{k})  \right)  \Bigr) , 
\EQNY
\BQNY
	 S_{s,t}^{\left(3\right) } 
	 \sim ME \left( 2\beta_{k}, \underline{\Pi} \left( \Psi \left( M_{s,t} (
	 \underline{Q}_{1}) \right), \ldots,\Psi \left( M_{s,t} (
	 \underline{Q}_{k}) \right)  \right) \right)  , \ 
   S_{s}^{\left( 2\right)}=S_{s,s}^{\left(3\right) },
\EQNY
and 
\BQNY 
M_{s,t} \left( \underline{Q}_{l}  \right) =
 		\left\{
			 \begin{array}{lcl}
			 V( \underline{Q}_{l} ) & \mbox{if} & l \in \left\lbrace s,t\right\rbrace  \notag \\
         		 \underline{Q}_{l}  & \mbox{otherwise} &          			 
              \end{array}              	
      \right. .
\EQNY	
\ES

\begin{korr}
Under the assumptions of Proposition \ref{prop1} it follows that $ S\sim ME(2\beta_k,\underline{P}), $ where the components of the vector of mixing weights $\underline{P}=(p_1,p_2,\ldots)$ are defined by
\BQNY
p_i&=&\xi_k \pi_i \left( \Psi(\underline{Q}_1),\ldots,\Psi(\underline{Q}_k) \right) 
  - \sum_{1 \leq s<t \leq k} \alpha_{s,t} \gamma_{s} \gamma_{t} \left[ 
\pi_i \left( \Psi \left( M_{s} ( \underline{Q}_{1}) \right), \ldots,\Psi \left( M_{s} (
	 \underline{Q}_{k}) \right)  \right) \right. 
		 \notag \\
&& \left. + \pi_i \left( \Psi \left( M_{t} ( \underline{Q}_{1}) \right), \ldots,\Psi \left( M_{t} (
	 \underline{Q}_{k}) \right)  \right)
- \pi_i \left( \Psi \left( M_{s,t} ( \underline{Q}_{1}) \right), \ldots,\Psi \left( M_{s,t} (
	 \underline{Q}_{k}) \right)  \right)\right] ,
\EQNY
where $ M_s=M_{s,s} $ and $ \pi_i $ are the components of $\underline{\Pi}$ defined in Remark \ref{rem:aggregME}.
\end{korr}
 
\begin{example} \label{example3.2} Bivariate mixed Erlang risks joined by Sarmanov's distribution.\\
Let $(X_1,X_2) \sim SME_2\left( \vk{\beta}=(\beta_1,\beta_2), (\underline{Q}_1,\underline{Q}_2),\alpha_{1,2}\right) $ with $\beta_1 < \beta_2$. It follows that $S= X_1+X_2 \sim ME(2\beta_2,\underline{P})$, where the components of the vector $\underline{P}$ are given below 
\BQNY
p_i&=&\left( 1+\alpha_{1,2} \gamma_1 \gamma_2\right) \pi_i\left( \Psi(\underline{Q}_1),\Psi(\underline{Q}_2)\right) 
- \alpha_{1,2} \gamma_1 \gamma_2 \left[ \pi_i\left( \Psi( V(\underline{Q}_1)), \Psi(\underline{Q}_2)\right) \right. \\
&&	\left. + \pi_i\left( \Psi( \underline{Q}_1), \Psi( V(\underline{Q}_2))\right) 
		-  \pi_i\left( \Psi( V(\underline{Q}_1)), \Psi( V(\underline{Q}_2))\right) \right] ,
\EQNY
such that $\sum_{i=1}^{\infty} p_i =1$.\\
\textbf{Numerical illustration.} As a numerical illustration, we considered a bivariate vector  $(X_1,X_2)$ such that 
$$ (X_1,X_2) \sim SME_2 \Bigl( \vk{\beta}= (0.9,0.95), \underline{Q}_1=(0.4,0.6), \underline{Q}_2=(0.8,0.2), \alpha_{1,2}=2.5\Bigr).$$
Thus, the densities of $X_1$ and $X_2$ can be, respectively, written as follows:
$$f_1(x)= 0.4 w_1(x,0.9)+0.6 w_2(x,0.9),\ ~ f_2(x)= 0.8 w_1(x,0.95)+0.2 w_2(x,0.95).$$
Moreover, from formula \eqref{eq: mean} we have
\BQNY
\mathbb{E} \left( f_i(X_i)\right) = \beta_i \sum_{l=1}^{2}  \sum_{j=1}^{2}   
	\begin{pmatrix} {l+j-2} \\ {l-1}  \end{pmatrix} \frac{q_{i,l}  q_{i,j}}{2^{l+j-1} },\ ~ i=1,2,
\EQNY
yielding $ \gamma_1=0.261,\gamma_2=0.3895$. In the following, we restrict to only two decimal places. Then the joint pdf of $ (X_1,X_2)$ is given by
$$ h(x_1,x_2)= f_1(x_1)f_2(x_2)\left( 1.25+2.5f_1(x_1)f_2(x_2)-0.97f_1(x_1)-0.65f_2(x_2)\right) .$$ 
 
Table \ref{table:quantmeas0} summarizes some quantitative measures related to the marginals $X_1$ and $X_2$.
\begin{center}
	\begin{tabular}{|c|c|c|c|c|}
\hline
& Expected value & Variance & Skewness & Kurtosis\\		
\hline
	$X_1$ & 1.78 & 2.27& 1.55& 6.50 \\
\hline
	$X_2$& 1.26 &1.51 &1.88 & 8.16\\
\hline
		\end{tabular}
	\captionof{table}{Quantitative measures for $X_1$ and $X_2$ (Example \ref{example3.2}).} \label{table:quantmeas0}
\end{center}
As stated above, the distribution of the aggregate risk $S$ is again mixed Erlang with scale parameter $2\beta_2=1.9$ and the mixing probabilities given in Table \ref{table:quantmeas1}.
\begin{center}
	\begin{tabular}{|c|c||c|c||c|c||c|c|}
\hline
 i &$p_i$ &i&$ p_i$&i& $p_i$ &i &$ p_i$ \\		
\hline
	1	&	0.0000   	&	11	&	 0.0262 	&	21	&	 0.0002 	&	31	&	8.635E-07\\
\hline
	2	&	 0.0827 	&	12	&	 0.0173 	&	22	&	 0.0001 	&	32	&	4.873E-07\\
\hline
	3	&	 0.1547 	&	13	&	 0.0112 	&	23	&	7.443E-05	&	 33 	&	2.743E-07\\
\hline
	4	&	 0.1709 	&	14	&	 0.0071 	&	24	&	4.326E-05	&	 34 	&	1.540E-07\\
\hline
	5	&	 0.1390 	&	15	&	 0.0045 	&	25	&	2.502E-05	&	 35 	&	8.625E-08\\
\hline
	6	&	 0.1162 	&	16	&	 0.0028 	&	26	&	1.441E-05	&	 36 	&	4.821E-08\\
\hline
	7	&	 0.0956 	&	17	&	 0.0017 	&	27	&	8.263E-06	&	 37 	&	2.689E-08\\
\hline
	8	&	 0.0744 	&	18	&	 0.0010 	&	28	&	4.722E-06	&	 38 	&	1.497E-08\\
\hline
	9	&	 0.0547 	&	19	&	 0.0006 	&	29	&	2.689E-06	&	 39 	&	8.319E-09\\
\hline
	10	&	 0.0385 	&	20	&	 0.0004 	&	30	&	1.526E-06	&	 40 	&	4.615E-09\\
\hline

		\end{tabular}
	\captionof{table}{Mixing probabilities of $S$ (Example \ref{example3.2}).} \label{table:quantmeas1}
\end{center}
\end{example}

We are now interested in quantifying the amount of capital $C_j(p)$ to be allocated to each risk $X_j,j \in I_k$.  

\BS\label{prop2}
Let $\textbf{X} \sim SME_k(\vk{\beta},\vk{\underline{Q}},\vk{\alpha})$ with $\beta_j \leq \beta_k,j\in I_k$, and let $s_p=VaR_p(S)$. Then the amount of capital $C_j$ allocated to each risk $X_j$ under the TVaR allocation principle as defined in \eqref{eq:ReinsTVar}  can be expressed as 
\BQN \label{eq:CTE}
C_j(p)=\frac{1}{1-p}\sum_{i=1}^{\infty}z_{i,j}\overline{W}_i(s_p,2\beta_k),
\EQN

where the mixing coefficients $z_{i,j}$ are given by (here the transform $\Psi$ is needed to obtain the common scale parameter $2\beta_k$)
\BQN \label{eq: Weight_Sn}
z_{i,j}&=& \xi_k \mu_j \pi_i \left(   \Psi\left( \widetilde{M}_j(\underline{Q}_1)\right), \ldots,  \Psi\left( \widetilde{M}_j(\underline{Q}_k)\right)  \right) -\sum_{1\leq a<b\leq k}\alpha_{a,b}\gamma_{a}\gamma_{b} \left[\varphi_{j;b} \pi_i \left( \Psi\left(\widetilde{M}_{j;b}(\underline{Q}_1)\right) , \ldots,  \Psi\left( \widetilde{M}_{j;b}(\underline{Q}_k)\right)  \right) \right. 
 \notag \\
&&\left.  +\varphi_{j;a} \pi_i \left( \Psi\left(\widetilde{M}_{j;a}(\underline{Q}_1)\right) , \ldots,  \Psi\left( \widetilde{M}_{j;a}(\underline{Q}_k)\right)  \right)
	-\varphi_{j;a,b} \pi_i \left( \Psi\left(\widetilde{M}_{j;a,b}(\underline{Q}_1)\right) , \ldots,  \Psi\left( \widetilde{M}_{j;a,b}(\underline{Q}_k)\right)  \right)	\right] ,				
\EQN
with $ \xi_k=1+ \sum_{1\leq a<b\leq k}\alpha_{a,b}\gamma_{a}\gamma_{b},\ \mu_i =\mathbb{E}( X_{i}) = \frac{1}{\beta_i} \sum_{k=1}^{\infty} kq_{i,k},\ \tilde{\mu}_i=\frac{1}{2\beta_i}\sum_{k=1}^{\infty} k v_{i,k}$ as defined in formula \eqref{mu_i}, 
\BQN \label{eq:varphis}
\varphi_{j;a,b}=\left\lbrace  
                \begin{array}{lcl}
         			 \mu_j & \mbox{if} & j\notin \{ a,b\}\\
         			 \tilde{\mu}_j & \mbox{if} & j \in \{a,b \}
              	\end{array},
      \right.
 \varphi_{j;a}=  \varphi_{j;a,a},  
\EQN

and 
\BQN \label{eq:Mtilde}
\widetilde{M}_{j}(\underline{Q}_i)=
\begin{cases}
\underline{Q}_i & \text{if }  i \neq j\\
G(\underline{Q}_i) & \text{if }  i=j\\
\end{cases}
, \widetilde{M}_{j;a,b}(\underline{Q}_i)=
\begin{cases}
\underline{Q}_i & \text{if }  i \notin \{j,a,b\}\\
V(\underline{Q}_i) & \text{if }  i=a \text{ and } i \notin \{j,b\} \text{ or} \\
 & \text{if } i= b \text{ and } i \notin \{j,a\} \\
G(\underline{Q}_i) & \text{if }  i=j \text{ and } i \notin \{a,b\}\\
G(V(\underline{Q}_i)) & \text{if } i= j=a \text{ and } i\neq b \text{ or}\\
 & \text{if } i= j=b \text{ and } i\neq a
\end{cases},
\widetilde{M}_{j;a}=\widetilde{M}_{j;a,a}.
\EQN
\ES

\begin{example} \label{example3.11} Capital allocation for bivariate mixed Erlang risks joined by Sarmanov's distribution.\\
In the bivariate case, with the above notation, $S_2=X_1+X_2$ is the aggregate risk of the portfolio and we consider $TVaR_p(S_2)$ to be the total capital needed to cover it, whereas $C_i$ is the part of this capital allocated to cover $X_i,i=1,2$. For a numerical illustration, we consider the bivariate vector used in Example \ref{example3.2}, but this time we vary the value of $\alpha_{1,2}$. Table \ref{table:quantmeas2} summarizes the results under the TVaR capital allocation principle assuming a tolerance level $p=99 \%$ (the second column shows the variance of $S_2$ denoted $\sigma^2_{S_2}$). \\

\begin{center}
	\begin{tabular}{|c|c|c|c|c|}
\hline
 $\alpha_{12}$&  $\sigma^2_{S_2}$&$C_1(99 \%)$ &$C_2(99 \%)$&$ TVaR_{99\%}(S_2) $\\	
\hline
3.4	& 4.0509&	 6.3920 	 & 	  4.3958	 & 	  10.7878	 \\ 
\hline
2.5	& 3.9788&	 6.3703 	 & 	 4.3556 	 & 	 10.7259 	 \\ 
\hline
1.5	& 3.8987&	 6.3458 	 & 	 4.3086 	 & 	10.6544	 \\ 
\hline
0.5	& 3.8186&	 6.3209 	 & 	 4.2589 	 & 	10.5798	 \\ 
\hline
0	& 3.7785&	6.3083	 & 	4.2330 	 & 10.5413		 \\ 
\hline
-0.5	& 3.7385&	6.2956 	 & 	4.2063  	 & 	10.5019	 \\ 
\hline
-1.5	& 3.6584&	 6.2698	 & 4.1505	  	 & 	10.4203	 \\ 
\hline
-2.1	& 3.6103&	6.2542  	 & 4.1154	  	 & 	10.3696	 \\ 
\hline
		\end{tabular}
	\captionof{table}{$TVaR_{99\%}(S_2)$ and capital allocated to each risk $X_i,i=1,2$ (Example \ref{example3.11}).} \label{table:quantmeas2}
\end{center}

It can be seen from Table \ref{table:quantmeas2} that the total capital needed to cover $S_2$  is dependent on $\alpha_{1,2}$ .  Actually, a larger  $\alpha_{1,2}$ implies a riskier portfolio (see the corresponding variance, $\sigma^2_{S_2}$) and thus, more capital is needed to cover each risk. Also, it can be seen that $X_1$ accounts for a larger capital than $X_2$ as it is riskier (having larger variance and expected value, see Table \ref{table:quantmeas0}).
\end{example}

\begin{appendices}
\appendix
\section{Dependence structure }

In this section, we discuss the dependence structure of two mixed Erlang distributed r.v.s $(X_1,X_2)$ joined  by the Sarmanov distribution with different kernel functions (in the insurance context, $X_1,X_2$ are dependent insurance risks). As before, the kernel functions are written in the form
$ \phi(x)= g(x)-\mathbb{E}(g(X)),$ with $g$ properly chosen. To model the dependence between the two r.v.s $X_1$ and $X_2$, we shall use Pearson's correlation coefficient denoted by $\rho_{1,2}$ and defined by 
\BQN\notag
\rho_{1,2}(X_1,X_2)= \frac{ \mathbb{E} (X_1 X_2) -\mathbb{E}(X_1)\mathbb{E}(X_2)}{\sigma_1 \sigma_2},
\EQN
where $\sigma_i=\sqrt{Var(X_i)},i=1,2$.
In the case of Sarmanov's distribution, $\rho_{1,2}$ can be rewritten as 
\BQN\label{eq:correlation}
\rho_{1,2}(X_1,X_2)= \frac{\alpha_{1,2} \mathbb{E}(X_1\phi_1(X_1))\mathbb{E}(X_2 \phi_2(X_2))}{\sigma_1 \sigma_2}.
\EQN
Based on (\ref{eq:correlation}), we hereafter present Pearson's correlation coefficient for different kernel functions along with its maximal and minimal values,  in the particular case of mixed Erlang marginals.\\ 

\textbf{Case 1:} Let $g(x)=f(x)$ (i.e., the marginal pdf), which leads to the kernel function $\phi(x)=f(x)-\mathbb{E}(f(X))$. Under the assumption $ X_i \sim ME(\beta_i,\underline{Q}_i),i=1,2,$ using Lemmas \ref{lem:tranformation ME} and \ref{lem:tranformation ME g=x}, we obtain
\BQNY
\mathbb{E}(X_i f_i(X_i))&=& \int_0^{\infty} x f(x,\beta_i,\underline{Q}_i)^2 \mathrm{d}x
= \gamma_i \int_0^{\infty} x \frac{f(x,\beta_i,\underline{Q}_i)^2}{\gamma_i} \mathrm{d}x 
= \gamma_i \int_0^{\infty} x f(x,2\beta_i,V(\underline{Q}_i)) \mathrm{d}x
= \gamma_i \widetilde{\mu}_i,
\EQNY
hence $\mathbb{E}(X_i\phi_i(X_i))=\gamma_i \widetilde{\mu}_i-\gamma_i \mu_i$ and Pearson's coefficient is now given by 
\BQN
\rho_{1,2}(X_1,X_2)= \frac{\alpha_{1,2} \gamma_1 \gamma_2(\widetilde{\mu}_1 -\mu_1)(\widetilde{\mu}_2 - \mu_2)}{\sigma_1 \sigma_2}.
\EQN
Thus, from \eqref{eq:Alpha12Bound}, the maximal and the minimal values of Pearson's correlation are respectively given by
\BQN
\rho_{1,2}^{max}(X_1,X_2)=
	   \frac{\gamma_1 \gamma_2(\widetilde{\mu}_1 -\mu_1)(\widetilde{\mu}_2 - \mu_2)}{\max \{ \gamma_{1} (M_2- \gamma_{2}),(M_1-\gamma_{1}) \gamma_{2} \} \sigma_1 \sigma_2},
\EQN
\BQN
\rho_{1,2}^{min}(X_1,X_2)=\frac{-\gamma_1 \gamma_2(\widetilde{\mu}_1 -\mu_1)(\widetilde{\mu}_2 - \mu_2)}{\max \{ \gamma_{1}  \gamma_{2},(M_1-\gamma_{1})(M_2- \gamma_{2}) \}\sigma_1 \sigma_2},
\EQN
where we recall $M_i=\underset{x \in \mathbb{R}}{\text{max }}f_i(x),i=1,2.$
\\

\textbf{Case 2:} We consider $g(x)=e^{-tx}$ for $t=1$, hence the corresponding kernel function is given by $\phi(x)= e^{-x} - \mathbb{E}\left( e^{-X}\right) $. Pearson's correlation coefficient along with its lower and upper bounds can be found in \cite{Hashorva_Rija14}, where the particular case of mixed Erlang marginals is emphasized.\\ 

\textbf{Case 3:} Let $g(x)=x^t$, in which case the kernel function is given by
$ \phi(x)= x^t - \mathbb{E}(X^t).$ A usual choice here is $t=1,$ which leads to the kernel $ \phi(x)= x - \mathbb{E}(X)$ and to the correlation $\rho_{1,2}(X_1,X_2)=\alpha_{1,2}\sigma_1 \sigma_2$. In this case, to fulfill the condition $1+\alpha_{1,2} \phi_{1}(x_{1})\phi_{2}(x_{2}) \geq 0$, upper truncated distributions can be considered for $X_1,X_2$, which is not the object of our study. However, if we denote by $T_i,i=1,2,$ the corresponding upper truncation points and consider that the marginal pdf's are defined only for non-negative values, then the maximal and the minimal values of the correlation coefficient are, respectively, given by
\BQNY
\rho_{1,2}^{max}(X_1,X_2)=
	   \frac{\sigma_1 \sigma_2}{\max \{ \mu_{1} (T_2- \mu_{2}),(T_1-\mu_{1}) \mu_{2} \}},\ 
\rho_{1,2}^{min}(X_1,X_2)=\frac{-\sigma_1 \sigma_2}{\max \{ \mu_{1}  \mu_{2},(T_1-\mu_{1})(T_2- \mu_{2}) \}}.
\EQNY \\

\textbf{Case 4:} We consider the FGM distribution already studied in \cite{Cossette_al13}, obtained for $g(x)=2\overline{F}(x)$, with the corresponding kernel function $\phi(x)=1-2F(x)$. Its Pearson's correlation coefficient is given by $ \rho_{1,2}=\frac{1}{3} \alpha_{1,2}.$

The minimal and maximal values of $\rho_{1,2}$ are  $-\frac{1}{3}$ and $\frac{1}{3}$, respectively, which is an important drawback of the FGM distribution. Moreover, in the particular case of mixed Erlang marginals joined by the FGM distribution, the Pearson correlation coefficient can be found in \cite{Cossette_al13}.

\begin{example}\label{exampleA} Comparison of the dependency  between two mixed Erlang r.v.s joined by the Sarmanov distribution.\\
a) We consider the bivariate random vector $$ (X_1,X_2) \sim SME_2 \Bigl( \vk{\beta}= (2,2.5), \underline{Q}_1=(0.45,0.55), \underline{Q}_2=(0.5,0.5), \alpha_{1,2}\Bigr).$$ In this application, we would like to compare the dependency between $X_1$ and $X_2$ based on the four different kernel functions described above. Therefore, we compute the upper and lower bounds of the  Pearson correlation coefficients for each kernel, together with the corresponding parameter $\alpha$, as summarized in the table below. As dicussed above, in Case 3 we considered upper truncated distributions with $T_1=T_2=15$ such that the tail functions of $ X_1,X_2$ in this truncation point are very small, hence making this case comparable with the other not-truncated ones.
It can be seen that the largest range of dependence corresponds to the kernel considered in \textbf{Case 1} (and studied in this paper) and the smallest to the truncated \textbf{Case 3}.\\

\begin{center}
	\begin{tabular}{|c|c|c|c|c|c|c|}
\hline
 &Kernel  &  $\alpha_{max} $&$\rho_{max}$ &$\alpha_{min}$&$ \rho_{min} $\\	
\hline
\textbf{Case 1}&$f(x)-\mathbb{E}(f(X))$	&	3.2100	 &  0.3023	 &-2.1289	  & -0.2005	 \\ 
\hline
\textbf{Case 2}&$e^{-x}- \mathbb{E}(e^{-X})$	&	3.5854	 & 0.1921		 & -3.000	   & -0.1607	 \\ 
\hline
\textbf{Case 3} &$x - \mathbb{E}(X)$	&	0.0896	 & 	0.0318	 & 	-0.0049 &  -0.0017  \\ 
\hline
\textbf{Case 4} &$1-2{F}(x)$	& 1.0000	 	 &0.2711	 & 	-1.0000 &-0.2711	 \\ 
\hline
		\end{tabular}
	\captionof{table}{Upper and lower bounds of Pearson's correlation coefficient for different kernel functions (Example \ref{exampleA}.a).} \label{table:quantmeas3}
\end{center}

b) In the sequel, we assumed a common scale parameter $\beta$ for both marginals and we plotted the upper and lower bounds of the correlation coefficient as a function of $\beta$ for the four kernel functions.\\
The figure below shows that in Case 1, the dependency increases with $\beta$, in contrast with Case 3 (considered with $T_1=T_2=15$), where the dependency decreases with $\beta$ quite rapidly from the maximum correlation coefficient to approximately 0. Case 2 and Case 4 show an almost constant dependency structure with respect to $\beta$.

\begin{figure} [H]
	\centering
	\includegraphics [scale=0.55]{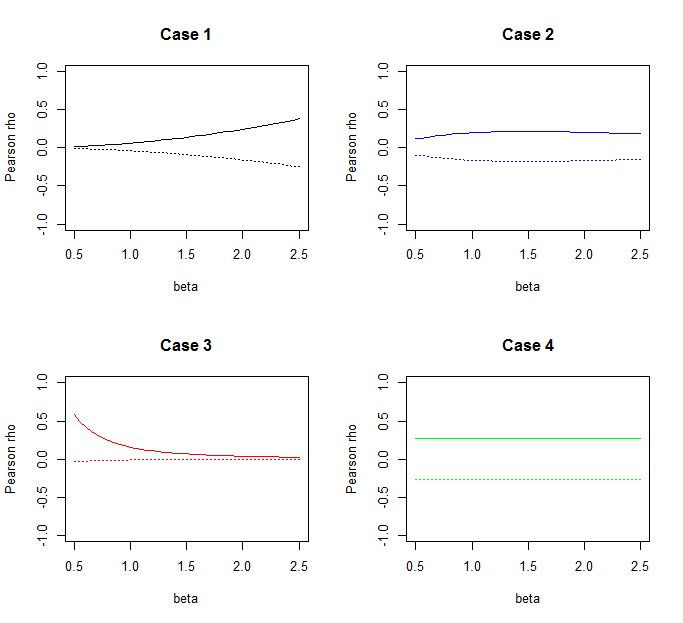}
	\caption{Peason's correlation coefficient for the same $\beta$ and different kernel functions (Example \ref{exampleA}.b). }
		\label{Figure: Premium_range}
	\end{figure} 
\end{example}


\section{Proofs} \label{SectionProofs}
 \textbf{Proof of Lemma} \ref{lem:tranformation ME}
We have 
\BQN \label{eq: pdf2}
f(x,\beta,\underline{Q})^2 &=& \left( \sum_{i=1}^{\infty} q_i w_i(x,\beta) \right)  \left( \sum_{j=1}^{\infty} q_j w_j (x,\beta) \right)  \notag \\
&=& \sum_{i=1}^{\infty}\sum_{j=1}^{\infty}q_i  q_j \frac{\beta^{i+j} x^{i+j-2} e^{-2\beta x}}{(i-1)!(j-1)!} \notag\\
&=& \sum_{i=1}^{\infty}\sum_{k=i}^{\infty}q_i  q_{k+1-i} \frac{\beta^{k+1} x^{k-1} e^{-2\beta x}}{(i-1)!(k-i)!} \notag\\
&=&\beta \sum_{k=1}^{\infty}\sum_{i=1}^{k}\begin{pmatrix} {k-1} \\ {i-1} \end{pmatrix}\frac{q_i  q_{k+1-i}}{2^{k} } w_{k}(x,2\beta).
\EQN

Also,
\BQN \label{eq: mean}
\mathbb{E} \left( f(X,\beta, \underline{Q})\right) &=&\int_{0}^{\infty} f(x,\beta,\underline{Q})f(x,\beta,\underline{Q}) \mathrm{d}x \notag \\ 
&=& \sum_{i=1}^{\infty}  \sum_{j=1}^{\infty} q_i  q_j \frac{\beta^{i+j}}{(i-1)!(j-1)!} \int_{0}^{\infty}  x^{i+j-2} e^{-2\beta x}\mathrm{d}x \notag \\ 
&=& \sum_{i=1}^{\infty}  \sum_{j=1}^{\infty} q_i  q_j \frac{\beta^{i+j}}{(i-1)!(j-1)!} \frac{(i+j-2)!}{(2\beta)^{i+j-1}} \notag \\  
&=& \beta \sum_{i=1}^{\infty}  \sum_{j=1}^{\infty}   
	\begin{pmatrix} {i+j-2} \\ {i-1}
 \end{pmatrix}\frac{q_i  q_j}{2^{i+j-1} }.
\EQN

Therefore, dividing \eqref{eq: pdf2} by \eqref{eq: mean}, we obtain
\BQNY
c(x, \beta, \underline{Q})&=& \frac{\beta \sum_{k=1}^{\infty}\sum_{i=1}^{k}\begin{pmatrix} {k-1} \\ {i-1} \end{pmatrix}\frac{q_i  q_{k+1-i}}{2^{k} } w_{k}(x,2\beta)}{\beta \sum_{i=1}^{\infty}  \sum_{j=1}^{\infty}   \begin{pmatrix} {i+j-2} \\ {i-1} \end{pmatrix}\frac{q_i  q_j}{2^{i+j-1} }}= \sum_{k=1}^{\infty} v_k w_k(x, 2\beta), 
\EQNY
where the coefficients $v_k$ are defined in \eqref{eq:vk}. \QED

\textbf{Proof of Lemma} \ref{lem:Aggreinsrisk}
We prove the result by induction. When $n=1$, from Lemma \ref{lem:Reinsriskdf} we have 
$F_{R_{1}}\left( y\right)=F_{Y_{1}}\left( y\right) =F_{X_{1}}\left( d_{1}\right) +H_{X_{1}}\left(
y,d_{1}\right) ,$ i.e., formula (\ref{FRn}) for $n=1$. Assuming now that the
formula (\ref{FRn}) holds for $n-1,$\ for $n$\ we obtain
\BQN
F_{R_{n}}\left( y\right) =\mathbb{P} \left( R_{n}\leq y\right) =\mathbb{P} \left(
R_{n}=0\right) +\mathbb{P} \left( 0<R_{n-1}+Y_{n}\leq y\right) .  \label{f0}
\EQN
But
$$ \mathbb{P} \left( R_{n}=0\right) =\mathbb{P} \left( Y_{i}=0,i\in I_{n}\right) =\mathbb{P} \left(
X_{i}\leq d_{i},i\in I_{n}\right) =\prod\limits_{i=1}^{n}F_{X_{i}}\left(
d_{i}\right) ,$$
while, using the induction hypothesis,
\BQN
\mathbb{P} \left( 0<R_{n-1}+Y_{n}\leq y\right)  &=&\mathbb{P} \left( Y_{n}=0,0<R_{n-1}\leq
y\right) +\mathbb{P} \left( R_{n-1}=0,0<Y_{n}\leq y\right)   \notag \\
&&+\mathbb{P} \left( 0<R_{n-1},0<Y_{n},0<R_{n-1}+Y_{n}\leq y\right)   \notag \\
&=&F_{X_{n}}\left( d_{n}\right) \sum_{k=1}^{n-1}\sum_{1\leq
j_{1}<...<j_{k}\leq n-1}H_{X_{j_{1}}+...+X_{j_{k}}}\left(
y,d_{j_{1}},...,d_{j_{k}}\right) \prod\limits_{i\in I_{n-1}\setminus \left\{
j_{1},...,j_{k}\right\} }F_{X_{i}}\left( d_{i}\right)   \notag \\
&&+\left( \prod\limits_{i=1}^{n-1}F_{X_{i}}\left( d_{i}\right)\right)  H_{X_{n}}\left(y,
d_{n}\right) +I_{R_{n}},  \label{f1}
\EQN
where, with $f_{Y_n}$ denoting the pdf of the r.v. $Y_n$,
\BQN
I_{R_{n}}&=&\mathbb{P} \left( 0<R_{n-1},0<Y_{n},0<R_{n-1}+Y_{n}\leq y\right)  \notag \\
 &=&\int_{0}^{y}\left[ \sum_{k=1}^{n-1}\sum_{1\leq
j_{1}<...<j_{k}\leq n-1}H_{X_{j_{1}}+...+X_{j_{k}}}\left(
y-u,d_{j_{1}},...,d_{j_{k}}\right) \prod\limits_{i\in I_{n-1}\setminus
\left\{ j_{1},...,j_{k}\right\} }F_{X_{i}}\left( d_{i}\right) \right]
f_{Y_{n}}\left( u\right) du  \notag \\
&=&\sum_{k=1}^{n-1}\sum_{1\leq j_{1}<...<j_{k}\leq n-1}\left(
\prod\limits_{i\in I_{n-1}\setminus \left\{ j_{1},...,j_{k}\right\}
}F_{X_{i}}\left( d_{i}\right) \right)
\int_{0}^{y}H_{X_{j_{1}}+...+X_{j_{k}}}\left(
y-u,d_{j_{1}},...,d_{j_{k}}\right) H'_{X_{n}}\left( u,d_{n}\right) du  \notag \\
&=&\sum_{k=1}^{n-1}\sum_{1\leq j_{1}<...<j_{k}\leq n-1}\left(
\prod\limits_{i\in I_{n-1}\setminus \left\{ j_{1},...,j_{k}\right\}
}F_{X_{i}}\left( d_{i}\right) \right)\sum_{h_{j_1}=0}^{\infty
}...\sum_{h_{j_k}=0}^{\infty }\Delta _{h_{j_1}}\left( d_{j_1},\beta ,\underline{Q}_{j_1}\right) 
\cdot ...\cdot \Delta _{h_{j_k}}\left( d_{j_k},\beta ,\underline{Q}_{j_k}\right) \notag \\
&\times & \int_{0}^{y} W_{\sum_{i=1}^{k}h_{j_i}+k}\left( y-u,\beta \right) 
\sum_{h_{n}=0}^{\infty }\Delta _{h_{n}}\left( d_{n},\beta ,\underline{Q}_{n}\right)
w_{h_n+1}(u,\beta)du
\notag \\
&=&\sum_{k=1}^{n-1}\sum_{1\leq j_{1}<...<j_{k}\leq n-1}\left(
\prod\limits_{i\in I_{n-1}\setminus \left\{ j_{1},...,j_{k}\right\}
}F_{X_{i}}\left( d_{i}\right) \right)
H_{X_{j_{1}}+...+X_{j_{k}}+X_{n}}\left(y,
d_{j_{1}},...,d_{j_{k}},d_{n}\right) .  \label{f2}
\EQN
For the last equality, apart the definition of $H$, we also used the fact that the convolution of two Erlang distributions having the same scale parameter is again an Erlang distribution with the same scale parameter, while its shape parameter equals the sum of the shape parameters of the convoluted distributions. Inserting now (\ref{f2})  into (\ref{f1}) and the result into (\ref{f0}) yields (\ref{FRn}). To obtain the formula of $ U ,$ we use
\BQNY
U_{\sum_{i=1}^{k}X_{i},X_{k+1}}(y,d_1,\ldots,d_{k+1})=\sum_{h_{1}=0}^{\infty
}...\sum_{h_{k+1}=0}^{\infty }\Delta _{h_{1}}\left( d_{1},\beta ,\underline{Q}_{1}\right) 
\cdot ...\cdot \Delta _{h_{k+1}}\left( d_{k+1},\beta ,\underline{Q}_{k+1}\right) J,
\EQNY
where
\BQNY
J&=&\int_y^ \infty \int_0^s u w_{h_{k+1}+1}(u,\beta)w_{\sum_{i=1}^{k}h_i+k}\left( s-u,\beta \right) duds \notag \\
&=& \frac{h_{k+1}+1}{\beta} \int_y^ \infty \int_0^s w_{h_{k+1}+2}(u,\beta)w_{\sum_{i=1}^{k}h_i+k}\left( s-u,\beta \right) duds\notag \\
&=& \frac{h_{k+1}+1}{\beta} \int_y^ \infty w_{\sum_{i=1}^{k+1}h_i+k+2}\left( s,\beta \right) ds,
\EQNY
which inserted into the above formula of $ U $ immediately yields the result.\QED

 \prooftheo{theorem:JointVectors}
The joint density of  $\vk{S}=(S_{1},\ldots, S_{n})$ is determined in terms of the joint density $h$ of $\left( X_1^{(1)},\ldots, X_{k_1}^{(1)};\ldots; X_{1}^{(n)},\ldots, X_{k_n}^{(n)} \right) $ as follows
\BQN \label{eq:defN}
f_{\vk{S}}(s_1,\ldots, s_n)= \idotsint\limits_{\left\{ \vk{x}=(x_1^{(1)},\ldots, x_{k_1}^{(1)} , \ldots, x_{1}^{(n)},\ldots,   x_{k_n}^{(n)}) \left\vert
\sum_{j=1}^{k_{i}} x_j^{(i)}=s_{i},i\in I_n \right. \right\} }  h(\x)d\vk{x}.
\EQN

Based on \eqref{eq: Sarmanovrv}, 
\BQNY  
h(\vk{x})
&=&\prod_{i=1}^{n}\prod_{j=1}^{k_{i}}f_{j}^{(i)}\left( x_{j}^{(i)}\right) +\sum_{1\leq
a<b\leq n}\sum_{s=1}^{k_{a}}\sum_{t=1}^{k_{b}}\alpha _{s,t}^{(a,b)}\left( \phi
_{s}^{(a)}f_{s}^{(a)}\right) \left( x_{s}^{(a)}\right) \left( \phi
_{t}^{(b)}f_{t}^{(b)}\right) \left( x_{t}^{(b)}\right) \underset{(i,j)\notin
\left\{ \left( a,s\right) ,\left( b,t\right) \right\} }{\prod%
\limits_{i=1}^{n}\prod\limits_{j=1}^{k_{i}}}f_{j}^{(i)}\left( x_{j}^{(i)}\right)  \\
&&+\sum_{a=1}^{n}\sum_{1\leq s<t\leq k_{a}}\alpha _{s,t}^{(a)}\left( \phi
_{s}^{(a)}f_{s}^{(a)}\right) \left( x_{s}^{(a)}\right) \left( \phi
_{t}^{(a)}f_{t}^{(a)}\right) \left( x_{t}^{(a)}\right) \underset{(i,j)\notin
\left\{ \left( a,s\right) ,\left( a,t\right) \right\} }{\prod%
\limits_{i=1}^{n}\prod\limits_{j=1}^{k_{i}}}f_{j}^{(i)}\left( x_{j}^{(i)}\right)  \\
&=&\prod_{i=1}^{n}\prod_{j=1}^{k_{i}}f_{j}^{(i)}\left( x_{j}^{(i)}\right) +\sum_{1\leq
a\leq b\leq n}\sum_{s=1}^{k_{a}}\sum_{t=T\left( s,a,b\right) }^{k_{b}}\alpha
_{s,t}^{(a,b)}\prod_{i=1}^{n}\prod_{j=1}^{k_{i}}\tilde{f}_{j;s,t}^{(i;a,b)}\left( x_{j}^{(i)}\right) ,
\EQNY
where $T\left( s,a,b\right) =\max \left\{ 1,\left( s+1) \mathbbm{1}_{\left( a=b\right)
}\right) \right\},\alpha _{s,t}^{(a,a)}=\alpha _{s,t}^{(a)}$ and,
for  $i\in I_n, j=1,\ldots,{k_i},$ 
$$\tilde{f}_{j;s,t}^{(i;a,b)}\left( x\right) =\left\{ 
\begin{array}{ll}
\left( \phi _{j}^{(i)}f_{j}^{(i)}\right) \left( x\right)  & if\ (i,j)\in
\left\{ \left( a,s\right) ,\left( b,t\right) \right\}  \\ 
f_{j}^{(i)}\left( x\right)  & otherwise.%
\end{array}%
\right. .$$
Therefore, we can express \eqref{eq:defN} as 
\BQNY
f_{\vk{S}}(\vk{s}) &=&\prod_{i=1}^{n}\idotsint\limits_{\mathbb{R}^{k_{i}-1}}\left( \prod_{j=1}^{k_{i}-1}f_{j}^{(i)}\left( x_{j}^{(i)}\right) \right)
f_{k_{i}}^{(i)}\left( s_{i}-\sum_{j=1}^{k_{i}-1}x_{j}^{(i)}\right)
dx_{1}^{\left( i\right) }...dx_{k_{i}-1}^{\left( i\right) } \\
&&+\sum_{1\leq a \leq b\leq n}\sum_{s=1}^{k_{a}}\sum_{t=T\left( s,a,b\right)
}^{k_{b}}\alpha _{s,t}^{(a,b)}\prod_{i=1}^{n}\idotsint\limits_{\mathbb{R}%
^{k_{i}-1}}\left( \prod_{j=1}^{k_{i}-1}\tilde{f}%
_{j;s,t}^{(i;a,b)}\left( x_{j}^{(i)}\right) \right) \tilde{f}_{k_{i};s,t}^{(i;a,b)}\left(
s_{i}-\sum_{j=1}^{k_{i}-1}x_{j}^{(i)}\right) dx_{1}^{\left( i\right)
}...dx_{k_{i}-1}^{\left( i\right) } \\
&=&\prod_{i=1}^{n}f_{S_{i}}(s_{i})+\sum_{1\leq a\leq b\leq
n}\sum_{s=1}^{k_{a}}\sum_{t=T\left( s,a,b\right) }^{k_{b}}\alpha _{s,t}^{(a,b)}\prod_{i=1}^{n}\tilde{f}_{S_{i;s,t}^{\left( a,b\right) }}(s_{i}),
\EQNY
with $ f_{S_{i}} $ and $ \tilde{f}_{S_{i;s,t}^{\left( a,b\right)}} $ defined in \eqref{eq:convol}. This completes the proof.
\QED\\

\proofprop{proposition:survival}	
The  df of $ \vk{S} $  is determined in terms of the joint pdf $h$ of  $\vk{X}$ as follows
\BQN \label{eq:defN1}
		F_{\vk{S}}(\vk{s})=	\mathbb{P}(S_1 \leqslant s_1, \ldots, S_n  \leqslant  s_n)=
	 		\idotsint\limits_{\left\{ \vk{x}=(x_1^{(1)},\ldots, x_{k_1}^{(1)} , \ldots, x_{1}^{(n)},\ldots,   x_{k_n}^{(n)}) \left\vert
\sum_{j=1}^{k_{i}} x_j^{(i)}\leq s_{i},i\in I_n \right. \right\} }  h(\x)d\vk{x}.
\EQN

Starting from \eqref{eq: Sarmanovrv}, the joint density of $\vk{X}$ is now given by 
\BQNY
h\left(\vk{x}\right)
&=&\prod\limits_{i=1}^{n}\prod\limits_{j=1}^{k_{i}}f_{j}^{\left( i\right)
}\left( x_{j}^{\left( i\right) }\right) \left[ 1+\sum_{1\leq a\leq b\leq
n}\sum_{s=1}^{k_{a}}\sum_{t=T\left( s,a,b\right) }^{k_{b}}\alpha _{s,t}^{(a,b)}\left(
f_{s}^{\left( a\right) }\left( x_{s}^{\left( a\right) }\right) -\gamma
_{s}^{\left( a\right) }\right) \left( f_{t}^{\left( b\right) }\left(
x_{t}^{\left( b\right) }\right) -\gamma _{t}^{\left( b\right) }\right) %
\right]  \\
&=&\prod\limits_{i=1}^{n}\prod\limits_{j=1}^{k_{i}}f_{j}^{\left( i\right)
}\left( x_{j}^{\left( i\right) }\right) \left( 1+\sum_{1\leq a\leq b\leq
n}\sum_{s=1}^{k_{a}}\sum_{t=T\left( s,a,b\right) }^{k_{b}}\alpha _{s,t}^{(a,b)}\gamma
_{s}^{\left( a\right) }\gamma _{t}^{\left( b\right) }\right)  \\
&&-\sum_{1\leq a\leq b\leq n}\sum_{s=1}^{k_{a}}\sum_{t=T\left( s,a,b\right)
}^{k_{b}}\alpha _{s,t}^{(a,b)}\gamma _{s}^{\left( a\right) }\gamma _{t}^{\left(
b\right) }\prod\limits_{i=1}^{n}\prod\limits_{j=1}^{k_{i}}\left(
f_{j;s}^{\left( i;a\right) }\left( x_{j}^{\left( i\right) }\right)
+f_{j;t}^{\left( i;b\right) }\left( x_{j}^{\left( i\right) }\right) \right) 
\\
&&+\sum_{1\leq a\leq b\leq n}\sum_{s=1}^{k_{a}}\sum_{t=T\left( s,a,b\right)
}^{k_{b}}\alpha _{s,t}^{(a,b)}\gamma _{s}^{\left( a\right) }\gamma _{t}^{\left(
b\right)
}\prod\limits_{i=1}^{n}\prod\limits_{j=1}^{k_{i}}f_{j;s,t}^{\left(
i;a,b\right) }\left( x_{j}^{\left( i\right) }\right) ,
\EQNY%
where, for $i\in I_n,j=1,...,k_{i},$%
\BQNY
f_{j;s,t}^{\left( i;a,b\right) }\left( x\right)  &=&\left\{ 
\begin{array}{ll}
\left. \left( f_{j}^{\left( i\right) }\left( x\right) \right) ^{2}\right/
\gamma _{j}^{\left( i\right) } & if\ \left( i,j\right) \in \left\{ \left(
a,s\right) ,\left( b,t\right) \right\}  \\ 
f_{j}^{\left( i\right) }\left( x\right)  & otherwise%
\end{array}%
\right. , \\
f_{j;s}^{\left( i;a\right) }\left( x\right)  &=&f_{j;s,s}^{\left(
i;a,a\right) }\left( x\right) .	
 \EQNY
By Lemma \ref{lem:tranformation ME}, $\frac {\left( f^{(i)}_j\right)^2}{\gamma^{(i)}_j}$ 
is the pdf of a mixed Erlang distribution (with twice the scale parameter), therefore, using also the notation $\xi_n$ from \eqref{eq:csi},
one can write \eqref{eq:defN1} as a sum-product of convolutions of mixed Erlang distributions as follows\\
\BQN \label{eq:conv}
		F_{\vk{S}}\left( \vk{s}\right)  &=&\xi_n
\prod\limits_{i=1}^{n}\int_{0}^{s_{i} }\int_{0}^{s_{i}-x_{1}^{\left(
i\right) } }...\int_{0}^{s_{i}-\sum_{j=1}^{k_{i}-2}x_{j}^{\left(
i\right) } }\prod\limits_{j=1}^{k_{i}-1}f_{j}^{\left( i\right)
}\left( x_{j}^{\left( i\right) }\right) F_{k_i}^{\left( i\right) }\left(
s_{i}-\sum_{j=1}^{k_{i}-1}x_{j}^{\left( i\right) }\right)
dx_{k_{i}-1}^{\left( i\right) }...dx_{1}^{\left( i\right) }\notag \\
&&-\sum_{1\leq
a\leq b\leq n}\sum_{s=1}^{k_{a}}\sum_{t=T\left( s,a,b\right) }^{k_{b}}\alpha _{s,t}^{(a,b)}\gamma _{s}^{\left( a\right) }\gamma _{t}^{\left( b\right) } \notag \\
&&\times \prod\limits_{i=1}^{n}\int_{0}^{s_{i} }\int_{0}^{s_{i}-x_{1}^{\left(
i\right) } }...\int_{0}^{s_{i}-\sum_{j=1}^{k_{i}-2}x_{j}^{\left(
i\right) } }\left[\prod\limits_{j=1}^{k_{i}-1} f_{j;s}^{\left( i;a\right) }\left(
x_{j}^{\left( i\right) }\right) F_{k_i;s}^{\left( i;a\right) }\left(
s_{i}-\sum_{j=1}^{k_{i}-1}x_{j}^{\left( i\right) }\right) \right. \notag \\
&&\left. +\prod\limits_{j=1}^{k_{i}-1}f_{j;t}^{\left( i;b\right) }\left( x_{j}^{\left( i\right) }\right)
F_{k_i;t}^{\left( i;b\right) }\left( s_{i}-\sum_{j=1}^{k_{i}-1}x_{j}^{\left(
i\right) }\right) \right] dx_{k_{i}-1}^{\left( i\right) }...dx_{1}^{\left(
i\right) }+\sum_{1\leq a\leq b\leq n}\sum_{s=1}^{k_{a}}\sum_{t=T\left(
s,a,b\right) }^{k_{b}}\alpha _{s,t}^{(a,b)}\gamma _{s}^{\left( a\right) }\gamma
_{t}^{\left( b\right) } \notag \\
&&\times \prod\limits_{i=1}^{n}\int_{0}^{s_{i} }\int_{0}^{s_{i}-x_{1}^{\left(
i\right) } }...\int_{0}^{s_{i}-\sum_{j=1}^{k_{i}-2}x_{j}^{\left(
i\right) } }\prod\limits_{j=1}^{k_{i}-1}f_{j;s,t}^{\left( i;a,b\right) }\left(
x_{j}^{\left( i\right) }\right) F_{k_i;s,t}^{\left( i;a,b\right) }\left(
s_{i}-\sum_{j=1}^{k_{i}-1}x_{j}^{\left( i\right) }\right)
dx_{k_{i}-1}^{\left( i\right) }...dx_{1}^{\left( i\right) }.
\EQN
Since $\beta_{k_n}^{(n)} \ge \beta_j^{(i)}, \forall j=1,\ldots, k_i,i\in I_n$, by Lemma \ref{lem:tranformation ME with different scale parameters} each $i$th mixed Erlang  component of \eqref{eq:conv}  can be transformed into a new mixed Erlang distribution  with a common scale parameter $2\beta_{k_n}^{(n)}$.
In addition, according to Remark \ref{rem:aggregME}, the convolution of mixed Erlang distributions belongs to the class of mixed Erlang distributions. Therefore, \eqref{eq:conv} can be expressed as a sum-product of mixed Erlang df's as follows
   \BQNY
		F_{\vk{S}}\left( s\right)=\xi_n
\prod\limits_{j=1}^{n}F_{S_{j}^{\left( 1 \right) }}\left(
s_{j}\right) -\sum_{1\leq a\leq b\leq n}\sum_{s=1}^{k_{a}}\sum_{t=T\left(
s,a,b\right) }^{k_{b}}\alpha _{s,t}^{(a,b)}\gamma _{s}^{\left( a\right) }\gamma
_{t}^{\left( b\right) }  \left( \prod\limits_{j=1}^{n}F_{S_{j;s}^{\left( 2;a\right)
}}\left( s_{j}\right) +\prod\limits_{j=1}^{n}F_{S_{j;t}^{\left(
2;b\right) }}\left( s_{j}\right)
-\prod\limits_{j=1}^{n}F_{S_{j;s,t}^{\left(3;a,b\right) }}\left(
s_{j}\right) \right) ,
	\EQNY
where $S_j^{(1)},S_{j;s}^{\left( 2;a\right)},S_{j;s,t}^{\left(3;a,b\right) }$ are defined by \eqref{eq:defS1}-\eqref{eq:defS2}. Thus the proof is complete. 
\QED

\proofprop{prop3}
The distribution of $R_{n}$ can be expressed in terms of the distribution of $\vk{S}$ as follows
\BQNY
F_{R_{n}}(y)&=&\mathbb{P}\left( \bigcap\limits_{i=1}^{n}\left( T_{i}=0\right) \right)
+\sum_{k=1}^{n}\sum_{1\leq j_{1}<...<j_{k}\leq n}\mathbb{P}\left(
\bigcap\limits_{i\in I_{n}\backslash \left\{ j_{1},...,j_{k}\right\} }\left(
T_{i}=0\right) \bigcap\limits_{i=1}^{k}\left( T_{j_{i}}>0\right) \bigcap
\left( \sum_{i=1}^{k}T_{j_{i}}\leq y\right) \right)  \\
&=&\mathbb{P}\left( \bigcap\limits_{i=1}^{n}\left( S_{i}\leq d_{i}\right)
\right) +\sum_{k=1}^{n}\sum_{1\leq j_{1}<...<j_{k}\leq n}\mathbb{P}\left(
\bigcap\limits_{i\in I_{n}\backslash \left\{ j_{1},...,j_{k}\right\} }\left(
S_{i}\leq d_{i}\right) \bigcap\limits_{i=1}^{k}\left(
S_{j_{i}}>d_{j_{i}}\right) \bigcap \left( \sum_{i=1}^{k}\left(
S_{j_{i}}-d_{j_{i}}\right) \leq y\right) \right). 
\EQNY
Considering now the df of $\vk{S}$ given in Proposition \ref{proposition:survival} and 
the definition of $H$ from Lemma \ref{lem:Aggreinsrisk}, we can rewrite $ F_{R_{n}} $ in the form \eqref{eq:FRnProp}, which completes the proof.
\QED

\proofprop{prop4}
For simplicity, we shall prove the case $ l=n $, the proof for a general $ l $ being similar, but with a notation more complicated. To use (\ref{eq:ReinsTVar}), we must evaluate
\BQN \label{eq:J12}
\mathbb{E}\left( T_{n}\mathbbm{1}_{\left\{ R_{n}>x_{p}\right\} }\right) =J_{1}+J_{2},
\EQN

where
\BQNY
J_{1} &=&\int_{x_{p}}^{\infty }uf_{T_{n},\left\{ R_{n-1}=0\right\} }\left(
u\right) du\text{ with }f_{T_{n},\left\{ R_{n-1}=0\right\} }\left( u\right) =%
\frac{\partial }{\partial u}\mathbb{P}\left( 0<T_{n}\leq u,R_{n-1}=0\right) ,
\\
J_{2} &=&\int_{x_{p}}^{\infty }\int_{0}^{s}uf_{T_{n},R_{n-1}}\left(
u,s-u\right) duds\text{ with }f_{T_{n},R_{n-1}}\left( u,v\right) =\frac{%
\partial ^{2}}{\partial u\partial v}\mathbb{P}\left( 0<T_{n}\leq
u,0<R_{n-1}\leq v\right) .
\EQNY

We shall now use Proposition \ref{prop3} and the notation from Lemma \ref{lem:Reinsriskdf}. We have
\BQNY
J_{1} &=&\int_{x_{p}}^{\infty }u\frac{\partial }{\partial u}\left[ \xi
_{n}\left( \prod\limits_{i=1}^{n-1}F_{S_{i}^{\left( 1\right) }}\left(
d_{i}\right) \right) H_{S_{n}^{\left( 1\right) }}\left( u,d_{n}\right)
-\sum_{1\leq a\leq b\leq n}\sum_{s=1}^{k_{a}}\sum_{t=T\left( s,a,b\right)
}^{k_{b}}\alpha _{s,t}^{\left( a,b\right) }\gamma _{s}^{\left( a\right)
}\gamma _{t}^{\left( b\right) }\right.  \\
&&\times \left. \left( H_{S_{n;s}^{\left( 2;a\right) }}\left( u,d_{n}\right)
\prod\limits_{i=1}^{n-1}F_{S_{i;s}^{\left( 2;a\right) }}\left( d_{i}\right)
+H_{S_{n;t}^{\left( 2;b\right) }}\left( u,d_{n}\right)
\prod\limits_{i=1}^{n-1}F_{S_{i;t}^{\left( 2;b\right) }}\left( d_{i}\right)
-H_{S_{n;s,t}^{\left( 3;a,b\right) }}\left( u,d_{n}\right)
\prod\limits_{i=1}^{n-1}F_{S_{i;s,t}^{\left( 3;a,b\right) }}\left(
d_{i}\right) \right)\right]  du \\
&=&\xi _{n}U_{S_{n}^{\left( 1\right) }}\left( x_{p},d_{n}\right)
\prod\limits_{i=1}^{n-1}F_{S_{i}^{\left( 1\right) }}\left( d_{i}\right)
-\sum_{1\leq a\leq b\leq n}\sum_{s=1}^{k_{a}}\sum_{t=T\left( s,a,b\right)
}^{k_{b}}\alpha _{s,t}^{\left( a,b\right) }\gamma _{s}^{\left( a\right)
}\gamma _{t}^{\left( b\right) } \\
&&\times \left( U_{S_{n;s}^{\left( 2;a\right) }}\left( x_{p},d_{n}\right)
\prod\limits_{i=1}^{n-1}F_{S_{i;s}^{\left( 2;a\right) }}\left( d_{i}\right)
+U_{S_{n;t}^{\left( 2;b\right) }}\left( x_{p},d_{n}\right)
\prod\limits_{i=1}^{n-1}F_{S_{i;t}^{\left( 2;b\right) }}\left( d_{i}\right)
-U_{S_{n;s,t}^{\left( 3;a,b\right) }}\left( x_{p},d_{n}\right)
\prod\limits_{i=1}^{n-1}F_{S_{i;s,t}^{\left( 3;a,b\right) }}\left(
d_{i}\right) \right) .
\EQNY

On the other hand, a reasoning similar with the one in the proof of Proposition \ref{prop3} yields
\BQNY
f_{T_{n},R_{n-1}}\left( u,v\right)  &=&\frac{\partial ^{2}}{\partial
u\partial v}\sum_{k=1}^{n-1}\sum_{1\leq j_{1}<...<j_{k}\leq n-1}\mathbb{P}%
\left( \left( 0<T_{n}\leq u\right) \bigcap\limits_{i\in I_{n-1}\backslash
\left\{ j_{1},...,j_{k}\right\} }\left( T_{i}=0\right)
\bigcap\limits_{i=1}^{k}\left( T_{j_{i}}>0\right) \bigcap \left(
\sum_{i=1}^{k}T_{j_{i}}\leq v\right) \right)  \\
&=&\frac{\partial ^{2}}{\partial u\partial v}\sum_{k=1}^{n-1}\sum_{1\leq
j_{1}<...<j_{k}\leq n-1}\left[ \xi _{n}H_{S_{n}^{\left( 1\right) }}\left(
u,d_{n}\right) H_{\sum_{i=1}^{k}S_{j_{i}}^{\left( 1\right) }}\left(
v,d_{j_{1}},..,d_{j_{k}}\right) \prod\limits_{i\in I_{n-1}\backslash \left\{
j_{1},...,j_{k}\right\} }F_{S_{i}^{\left( 1\right) }}\left( d_{i}\right)
\right.  \\
&&-\sum_{1\leq a\leq b\leq n}\sum_{s=1}^{k_{a}}\sum_{t=T\left( s,a,b\right)
}^{k_{b}}\alpha _{s,t}^{\left( a,b\right) }\gamma _{s}^{\left( a\right)
}\gamma _{t}^{\left( b\right) } \\
&&\times \left( H_{S_{n;s}^{\left( 2;a\right) }}\left( u,d_{n}\right)
H_{\sum_{i=1}^{k}S_{j_{i};s}^{\left( 2;a\right) }}\left(
v,d_{j_{1}},..,d_{j_{k}}\right) \prod\limits_{i\in I_{n-1}\backslash \left\{
j_{1},...,j_{k}\right\} }F_{S_{i;s}^{\left( 2;a\right) }}\left( d_{i}\right)
\right.  \\
&&+H_{S_{n;t}^{\left( 2;b\right) }}\left( u,d_{n}\right)
H_{\sum_{i=1}^{k}S_{j_{i};t}^{\left( 2;b\right) }}\left(
v,d_{j_{1}},..,d_{j_{k}}\right) \prod\limits_{i\in I_{n-1}\backslash \left\{
j_{1},...,j_{k}\right\} }F_{S_{i;t}^{\left( 2;b\right) }}\left( d_{i}\right) 
\\
&&\left. \left. -H_{S_{n;s,t}^{\left( 3;a,b\right) }}\left( u,d_{n}\right)
H_{\sum_{i=1}^{k}S_{j_{i};s,t}^{\left( 3;a,b\right) }}\left(
v,d_{j_{1}},..,d_{j_{k}}\right) \prod\limits_{i\in I_{n-1}\backslash \left\{
j_{1},...,j_{k}\right\} }F_{S_{i;s,t}^{\left( 3;a,b\right) }}\left(
d_{i}\right) \right) \right] ,
\EQNY
from where, using the notation introduced in Lemma \ref{lem:Aggreinsrisk}, we obtain
\BQNY
J_{2} &=&\sum_{k=1}^{n-1}\sum_{1\leq j_{1}<...<j_{k}\leq n-1}\left[ \xi
_{n}U_{\sum_{i=1}^{k}S_{j_{i}}^{\left( 1\right) },S_{n}^{\left( 1\right)
}}\left( x_{p},d_{j_{1}},..,d_{j_{k}},d_{n}\right) \prod\limits_{i\in
I_{n-1}\backslash \left\{ j_{1},...,j_{k}\right\} }F_{S_{i}^{\left( 1\right)
}}\left( d_{i}\right) \right.  \\
&&-\sum_{1\leq a\leq b\leq n}\sum_{s=1}^{k_{a}}\sum_{t=T\left( s,a,b\right)
}^{k_{b}}\alpha _{s,t}^{\left( a,b\right) }\gamma _{s}^{\left( a\right)
}\gamma _{t}^{\left( b\right) }\left( U_{\sum_{i=1}^{k}S_{j_{i};s}^{\left(
2;a\right) },S_{n;s}^{\left( 2;a\right) }}\left(
x_{p},d_{j_{1}},..,d_{j_{k}},d_{n}\right) \prod\limits_{i\in
I_{n-1}\backslash \left\{ j_{1},...,j_{k}\right\} }F_{S_{i;s}^{\left(
2;a\right) }}\left( d_{i}\right) \right.  \\
&&+U_{\sum_{i=1}^{k}S_{j_{i};t}^{\left( 2;b\right) },S_{n;t}^{\left(
2;b\right) }}\left( x_{p},d_{j_{1}},..,d_{j_{k}},d_{n}\right)
\prod\limits_{i\in I_{n-1}\backslash \left\{ j_{1},...,j_{k}\right\}
}F_{S_{i;t}^{\left( 2;b\right) }}\left( d_{i}\right)  \\
&&\left. \left. -U_{\sum_{i=1}^{k}S_{j_{i};s,t}^{\left( 3;a,b\right)
},S_{n;s,t}^{\left( 3;a,b\right) }}\left(
x_p,d_{j_{1}},..,d_{j_{k}},d_{n}\right) \prod\limits_{i\in I_{n-1}\backslash
\left\{ j_{1},...,j_{k}\right\} }F_{S_{i;s,t}^{\left( 3;a,b\right) }}\left(
d_{i}\right) \right) \right] .
\EQNY
Inserting now the formulas of $ J_1,J_2 $ into \eqref{eq:J12} and the result into \eqref{eq:ReinsTVar} yields the formula of $C_n$. Thus, the proof is complete.
\QED 
 
\proofprop{prop2}
Without loss of generality, we assume that $ j<k $. To prove the stated formula of $C_j$, we need to find 
\BQN\label{eq:expectedvalue}
\mathbb{E}\left( X_j\mathbbm{1}_{\{S>s_p\}}\right) =\int_{s_p}^\infty \int_0^s \int_0^{s-x_1}\ldots \int_0^{s-\sum_{i=1}^{k-2}x_{i}} x_j h\left( x_1,x_2,\ldots,x_{k-1},s-\sum_{i=1}^{k-1}x_{i}\right) \mathrm{d}{x_{k-1}}\ldots\mathrm{d}{x_2}\mathrm{d}{x_1}\mathrm{d} s,
\EQN

with $h$ as defined in \eqref{eq:Sarmanovdensity} with kernels \eqref{eq:condkern1}. We denote $ \xi_k=1+ \sum_{1\leq a<b\leq k}\alpha_{a,b}\gamma_{a}\gamma_{b}$ and first evaluate $x_j h(\x)$ as 
\BQNY 
x_j h(\x)&=& x_j\prod_{i=1}^{k}  f_{i}(x_i) \left[ 1+ \sum_{1\leq a<b\leq k}\alpha_{a,b}
 \left( f_a(x_a)f_b(x_b)-\gamma_{a}f_b(x_b)-\gamma_{b}f_a(x_a)+\gamma_{a}\gamma_{b}\right) \right] 
   \\ 
&=&\xi_k \mu_j \prod_{i=1}^{k}  f_{i,j}(x_i)-\sum_{1\leq a<b\leq k}\alpha_{a,b}\gamma_{a}\gamma_{b} \left[ \varphi_{j;b} \prod_{i=1}^{k}  f_{i,j;b}(x_i)+\varphi_{j;a} \prod_{i=1}^{k}  f_{i,j;a}(x_i)-
\varphi_{j;a,b}\prod_{i=1}^{k}  f_{i,j;a,b}(x_i)\right], 	
\EQNY
where $ \varphi_{j;a,b} $ and $\varphi_{j;a} $ are defined in \eqref{eq:varphis}, and, for $i \in I_k$, we define the following pdf's
\BQNY
f_{i,j}(x) &=&  \left\{
			 	   \begin{array}{lcl}
         			 f_i(x) & \mbox{if} & i \neq j \\
         			\frac { x f_{j} (x)}{\mu_j} & \mbox{if} & i = j
              	\end{array},
      \right.
	\\
f_{i,j;a,b}(x) &=&  \left\{
			 	   \begin{array}{lcl}
         			 f_i(x) & \mbox{if} & i \notin \{j,a,b\} \\
         			 \frac { x f_{j} (x)}{\mu_j} & \mbox{if} & i = j\notin \{a,b\}\\
         			\frac { f_i^2(x)}{\gamma_i} & \mbox{if} & i \in \{a,b\},i\neq j\\
                    \frac { x \left( f_{j}^2 (x)/ \gamma_j \right)}{\widetilde{\mu_j}}  & \mbox{if} & i = j \in \{a,b\} 	\end{array}  ,  \mbox{ while } f_{i,j;a}(x)=f_{i,j;a,a}(x).
      \right. 
 \EQNY
According to Lemmas \ref{lem:tranformation ME}, \ref{lem:tranformation ME g=x} and \ref{lem:tranformation ME with different scale parameters}, the just defined pdf's can be regarded of mixed Erlang type with parameter $2\beta_k$. Thus, (\ref{eq:expectedvalue}) becomes
 \BQNY 
\mathbb{E}\left( X_j\mathbbm{1}_{\{S>s_p\}}\right) &=&\xi_k \mu_j\int_{s_p}^\infty \int_0^s \int_0^{s-x_1}\ldots \int_0^{s-\sum_{i=1}^{k-2}x_{i}} f_{k,j}\left( s-\sum_{i=1}^{k-1}x_{i}\right) \prod_{i=1}^{k-1}  f_{i,j}(x_i) \mathrm{d}{x_{k-1}}\ldots\mathrm{d}{x_2}\mathrm{d}{x_1}\mathrm{d} s \\
&&-\sum_{1\leq a<b\leq k}\alpha_{a,b}\gamma_{a}\gamma_{b}\int_{s_p}^\infty \int_0^s \int_0^{s-x_1}\ldots \int_0^{s-\sum_{i=1}^{k-2}x_{i}} \left[
\varphi_{j;b} f_{k,j;b}\left( s-\sum_{i=1}^{k-1}x_{i}\right) \prod_{i=1}^{k-1} f_{i,j;b}(x_i) \right. \\
&& \left. +\varphi_{j;a} f_{k,j;a}\left( s-\sum_{i=1}^{k-1}x_{i}\right) \prod_{i=1}^{k-1}  f_{i,j;a}(x_i)-
\varphi_{j;a,b} f_{k,j;a,b}\left( s-\sum_{i=1}^{k-1}x_{i}\right) \prod_{i=1}^{k-1}  f_{i,j;a,b}(x_i)
\right] \mathrm{d}{x_{k-1}}\ldots\mathrm{d}{x_2}\mathrm{d}{x_1}\mathrm{d} s,
\EQNY
i.e., the sum of four integrals consisting of tails of convolutions of mixed Erlang distributions, which leads to the following four mixed Erlang distributions, respectively,
$$ME\left( 2\beta_k, \underline{\Pi} \left( \Psi\left( \widetilde{M}_j(\underline{Q}_1)\right) , \ldots,  \Psi\left( \widetilde{M}_j(\underline{Q}_k)\right)  \right)   \right), 
ME\left( 2\beta_k, \underline{\Pi} \left( \Psi\left( \widetilde{M}_{j;b}(\underline{Q}_1)\right) , \ldots,  \Psi\left( \widetilde{M}_{j;b}(\underline{Q}_k)\right)  \right)   \right),$$
$$ME\left( 2\beta_k, \underline{\Pi} \left( \Psi\left( \widetilde{M}_{j;a}(\underline{Q}_1)\right) , \ldots,  \Psi\left( \widetilde{M}_{j;a}(\underline{Q}_k)\right)  \right)   \right),
ME\left( 2\beta_k, \underline{\Pi} \left( \Psi\left( \widetilde{M}_{j;a,b}(\underline{Q}_1)\right) , \ldots,  \Psi\left( \widetilde{M}_{j;a,b}(\underline{Q}_k)\right)  \right)   \right),$$
where $\widetilde{M}_{j},\widetilde{M}_{j;a}$ and $\widetilde{M}_{j;a,b}$ are defined in \eqref{eq:Mtilde}.
Then formula \eqref{eq:CTE} holds with the mixing coefficients $z_{i,j}$ defined in (\ref{eq: Weight_Sn}). This completes the proof. \QED \\

\end {appendices}

\textbf{Acknowledgments}. 
G. Ratovomirija acknowledges partial support from the project RARE -318984
 (an FP7  Marie Curie IRSES Fellowship) and Vaudoise Assurances. 
 \bibliographystyle{plain}
\bibliography{aggSarmanovD}

\end{document}